\documentstyle[amssymb,12pt]{amsart}
\newtheorem{prop}{Proposition}[section]
\newtheorem{prop:def}{Proposition-Definition}[section]

\newtheorem{lemma}{Lemma}[section]
\newtheorem{thm}{Theorem}[section]

\theoremstyle{remark}
\newtheorem{remark}{Remark}

\begin{document} 

\newcommand{\nc}{\newcommand} \nc{\on}{\operatorname}

\nc{\¦}{{|}}

\nc{\pa}{\partial} 

\nc{\cA}{{\cal A}}\nc{\cB}{{\cal B}}\nc{\cC}{{\cal C}} \nc{\cE}{{\cal E}}
\nc{\cG}{{\cal G}}\nc{\cH}{{\cal H}}\nc{\cX}{{\cal X}}\nc{\cR}{{\cal R}}
\nc{\cL}{{\cal L}} \nc{\cK}{{\cal K}}\nc{\cO}{{\cal O}} 
\nc{\cF}{{\cal F}}\nc{\cM}{{\cal M}} \nc{\cW}{{\cal W}} 

\nc{\sh}{\on{sh}}\nc{\Id}{\on{Id}}\nc{\Diff}{\on{Diff}}
\nc{\ad}{\on{ad}}\nc{\Der}{\on{Der}}\nc{\End}{\on{End}}
\nc{\res}{\on{res}}\nc{\ddiv}{\on{div}} \nc{\FS}{\on{FS}}
\nc{\card}{\on{card}}\nc{\dimm}{\on{dim}}\nc{\gr}{\on{gr}}
\nc{\Jac}{\on{Jac}}\nc{\Ker}{\on{Ker}} \nc{\Den}{\on{Den}}
\nc{\Imm}{\on{Im}}\nc{\limm}{\on{lim}}\nc{\Ad}{\on{Ad}}
\nc{\ev}{\on{ev}} \nc{\Hol}{\on{Hol}}\nc{\Det}{\on{Det}}
\nc{\Cone}{\on{Cone}} \nc{\pseudo}{{\on{pseudo}}}
\nc{\class}{{\on{class}}}\nc{\rat}{{\on{rat}}}
\nc{\local}{{\on{local}}}\nc{\an}{{\on{an}}}
\nc{\Lift}{{\on{Lift}}}\nc{\Mer}{{\on{Mer}}}\nc{\mer}{{\on{mer}}}
\nc{\lift}{{\on{lift}}}\nc{\diff}{{\on{diff}}}\nc{\Aut}{{\on{Aut}}}
\nc{\DO}{{\on{DO}}}\nc{\Pic}{{\on{Pic}}}

\nc{\Bun}{\on{Bun}}\nc{\diag}{\on{diag}}\nc{\KZ}{{\on{KZ}}}
\nc{\CB}{{\on{CB}}}\nc{\out}{{\on{out}}}\nc{\Hom}{{\on{Hom}}}
\nc{\FO}{{\on{FO}}}

\nc{\al}{\alpha}\nc{\de}{\delta}\nc{\si}{\sigma}\nc{\ve}{\varepsilon}\nc{\z}{\zeta}
\nc{\vp}{\varphi} \nc{\la}{{\lambda}}\nc{\g}{\gamma}\nc{\eps}{\epsilon}
\nc{\PsiDO}{\Psi\on{DO}}

\nc{\AAA}{{\mathbb A}}\nc{\CC}{{\mathbb C}}\nc{\NN}{{\mathbb N}}
\nc{\PP}{{\mathbb P}}\nc{\RR}{{\mathbb R}}\nc{\VV}{{\mathbb V}}
\nc{\ZZ}{{\mathbb Z}}

\nc{\bla}{{\mathbf \lambda}} \nc{\kk}{{\mathbf k}} 
\nc{\bv}{{\mathbf v}}\nc{\bz}{{\mathbf z}}\nc{\bt}{{\mathbf t}} 
\nc{\bP}{{\mathbf P}}

\nc{\A}{{\mathfrak a}}\nc{\B}{{\mathfrak b}}\nc{\G}{{\mathfrak g}}
\nc{\HH}{{\mathfrak h}}\nc{\mm}{{\mathfrak m}}\nc{\N}{{\mathfrak n}}
\nc{\SG}{{\mathfrak S}}\nc{\La}{\Lambda}

\nc{\wt}{\widetilde}
\nc{\wh}{\widehat} \nc{\bn}{\begin{equation}}\nc{\en}{\end{equation}}
\nc{\SL}{{\mathfrak{sl}}}\nc{\ttt}{{\mathfrak{t}}}
\nc{\GL}{{\mathfrak{gl}}}

\title{Quantization of canonical cones of algebraic curves}

\begin{abstract}

We introduce a quantization of the graded algebra of functions 
on the canonical cone of an algebraic curve $C$, based on the theory of formal 
pseudodifferential operators. When $C$ is a complex curve with Poincar\'e 
uniformization, we propose another, equivalent construction, based on the 
work of Cohen-Manin-Zagier on Rankin-Cohen brackets. We give a presentation of 
the quantum algebra when $C$ is a rational curve, and discuss the problem of 
constructing algebraically "differential liftings".  

\end{abstract}

\author{B. Enriquez and A. Odesskii}

\address{B.E.: IRMA, Universit\'e  Louis Pasteur, 7, rue Ren\'e Descartes, 
67084 Strasbourg, France}

\address{A.O.: Landau Institute of Theoretical Physics, 2, Kosygina str., 
117334 Moscow, Russia}
\maketitle

\subsection*{Introduction}

To any pair $(C,D)$ of a curve and an effective divisor are associated
the morphism $C \to \PP(H^0(C,K(D))^*)$, where $K$ is the canonical 
bundle of $C$, and the corresponding cone $\Cone(C,D) \subset H^0(C,K(D))^*$.  
The function algebra of this cone is a graded algebra with Poisson structure. 
The purpose of this paper is to construct a quantization of this algebra.  
We will propose two equivalent solutions of this problem: 

(1) a solution based on the theory of formal pseudodifferential operators
(Section \ref{sect:pseudo}). Here the base field may be any algebraically 
closed field $\kk$ of characteristic zero. We show that the function algebras 
on $\Cone(C,D)$, as well as their quantizations, are functorial in 
the pair $(C,D)$ (Section \ref{sect:functoriality}). We also show (Section \ref{line}) 
that this construction can be twisted by a "generalized line bundle", i.e., 
an element of $\{$divisors with coefficients in $\kk\}/$ linear equivalence. 

(2) when the base field in $\CC$, we also present an analytic approach 
using Poincar\'e uniformization (Section \ref{sect:analytic}). This solution 
uses the results of \cite{CMZ} on Rankin-Cohen brackets (see also \cite{O}). 

In Section \ref{sect:example}, we give a presentation of the quantum 
algebra, when $C$ is a rational curve. 

In Section \ref{sect:lifting}, we discuss the problem of constructing
local, or differential, liftings from the classical algebra to the algebra of 
pseudodifferential operators. We show that Poincar\'e uniformization provides 
such liftings, in the analytic framework. We also discuss this problem in
the algebraic framework. 

In Section \ref{sect:final}, we discuss possible relations 
with the elliptic algebras of \cite{FO}, with Kontsevich 
quantization and with the problem of quantizing the Beauville 
hamiltonians of \cite{Beauville}. 

\section{Poisson algebras associated with canonical cones of curves}

\subsection{}

Let $C$ be a smooth, projective, connected complex curve (the constructions of 
this section can be generalized to the case where the base field is any algebraically 
closed field of characteristic zero). 
Let $K$ be its canonical  bundle. Let $D$ be an effective 
divisor of $C$; we set $D = \sum_{P\in C} \delta_P P$, where each 
$\delta_P$ is an integer $\geq 0$ and all but finitely many $\delta_P$ are zero. 
To these data is attached the morphism 
$$
C \to \PP(H^0(C,K(D))^*) 
$$ 
and the cone $\Cone(C,D)$, which is the preimage of $C$ by the map  
$H^0(C,K(D))^*\to \PP(H^0(C,K(D))^*)$. 
When $D = 0$, $\Cone(C,D)$ is the {\it canonical cone} of $C$. 
To each pair $D\geq D'$ is attached a morphism of cones 
$\Cone(C,D) \to \Cone(C,D')$. 

Moreover, the function ring of $\Cone(C,D)$ is a Poisson 
algebra. As an algebra, this is the graded ring 
\begin{equation} \label{graded:ring}
A^{(D)} = \bigoplus_{i\geq 0} H^0(C, K(D)^{\otimes i}); 
\end{equation}
we will denote by $A^{(D)}_i$ the part of $A^{(D)}$ of degree $i$.
For each $D$, we have an inclusion of graded rings $A^{(D)} \subset 
A^{\on{rat}}$, where
$$
A^{\on{rat}} = \bigoplus_{i\geq 0} \{\on{rational}\ i\on{-differentials}
\ \on{on}\ C\}. 
$$
We will define a Poisson structure on $A^{\on{rat}}$, which induces a 
Poisson structure on each $A^{(D)}$.  For this, we will choose a 
nonzero rational differential $\al$ on $C$, and define a Poisson structure 
$\{,\}_\al$; then we will show that this bracket is independent on the choice of $\al$. 

Let us denote by $\nabla^\al$ the meromorphic connection on $K^{\otimes i}$, 
such that if $\omega$ is a rational section of $K^{\otimes i}$, then 
$$
\nabla^\al(\omega) = \al^i d(\omega / \al^i). 
$$
Then we set, for $\omega,\omega'$ homogeneous of degrees $i,i'$,  
$$
\{\omega,\omega'\}_\al = i'\omega' \nabla^{\al}(\omega)
- i \omega \nabla^\al(\omega').  
$$ 

\begin{prop} \label{1.1}
The bracket $\{,\}_\al$ is independent on $\al$. We denote it 
by $\{,\}$. It is a Poisson bracket on $A^{\on{rat}}$, 
taking $A_i^{\on{rat}} \otimes A_j^{\on{rat}}$  to 
$A_{i+j+1}^{\on{rat}}$.  It restricts to a Poisson bracket 
on $A^{(D)}$. When the effective divisors $D_1$ and $D_2$ are linearly equivalent, 
the algebras $A^{(D_1)}$ and $A^{(D_2)}$ are isomorphic as graded 
algebras and as Poisson algebras.  
\end{prop}

{\em Proof.} Let us prove the independence on $\al$. Let $\beta$
be another differential. We have $\beta = F \al$, for $F$ a nonzero element of 
$\CC(C)$ (the field of rational functions on $C$). Then if $\omega$ 
is a rational section of $K^{\otimes i}$, we get 
$$
\nabla^{\beta}(\omega) = \nabla^{\al}(\omega) - i {{dF}\over F} \omega , 
$$  
so 
$$
\{\omega,\omega'\}_\beta = 
i'\omega' \nabla^{\beta}(\omega)
- i \omega \nabla^\beta(\omega')
= i'\omega' ( \nabla^{\al}(\omega)
- i {{dF}\over F} \omega) - \big((\omega,i) \leftrightarrow (\omega',i')\big) 
= \{\omega,\omega'\}_\al, 
$$
as $-ii'\omega\omega'{{dF}\over F}$ is symmetric under the exchange 
$(\omega,i) \leftrightarrow (\omega',i')$, so $\{\omega,\omega'\}_\beta = 
\{\omega,\omega'\}_\al$. We then define $\{\omega,\omega'\}$ as the common 
value of all $\{\omega,\omega'\}_\al$. 

It is easy to check that for any $\al$, $\{,\}_\al$ satisfies the Poisson 
bracket axioms, so the same is true for $\{,\}$. 

Let us show that $\{A_i^{(D)},A_{i'}^{(D)}\}\subset A^{(D)}_{i+i'+1}$. 
For this, we show that if $\omega$ (resp., $\omega'$) has a pole at 
$P$ of order $\leq i \delta_P$ (resp., $i'\delta_P$), then
$\{\omega,\omega'\}$ has a pole at $P$ of order 
$\leq \on{deg}(\{\omega,\omega'\})\delta_P = (i+i'+1)\delta_P$. 
Let $\al_P$ be a rational differential on $C$, such that 
$P$ is neither a zero nor a pole of $\al_P$. Then $\{\omega,\omega'\}
= \{\omega,\omega'\}_{\al_P}$. 
The terms of order $(i+i')\delta_P +1$ cancel each other, so the order of the 
pole of $\{\omega,\omega'\}_{\al_P}$ at $P$ is $\leq (i+i')\delta_P$. 
Since $(i+i')\delta_P\leq (i+i'+1)\delta_P$, we get 
$\{A_i^{(D)},A_{i'}^{(D)}\}\subset A^{(D)}_{i+i'+1}$. 

Finally, if $D_1 - D_2 = (f)$, where $f\in \CC(C)^\times$, then the $i$th 
component $A^{(D_1)}_i \to A^{(D_2)}_i$ of the isomorphism $A^{(D_1)} \to A^{(D_2)}$ 
takes $\omega \in H^0(C,K(D_1)^{\otimes i})$ to $\omega f^i \in 
H^0(C,K(D_2)^{\otimes i})$. 
\hfill \qed \medskip 

Then the natural morphism $A^{(D')} \hookrightarrow A^{(D)}$
attached to $D\geq D'$ is a Poisson algebra morphism, so 
$\Cone(C,D) \to \Cone(C,D')$ is a Poisson morphism. 

Moreover, one can describe the structure of symplectic leaves of 
$\Cone(C,D)$. Let us denote by $\on{Supp}(D)$ the support 
$\{P\in C\¦ \delta_P \neq 0\}$ of $D$.  

\begin{prop}
There exists a finite subset $D'$ of $C$, such that 
$\on{Supp}(D) \subset D' \subset \on{Supp}(D) \cup \{$Weierstrass points of 
$C\}$, with the following property. The symplectic leaves of $Cone(C,D)$ 
are of two types: 

-- each point of the preimage of $D'$ by $\Cone(C,D) \to C$ is a 
$0$-dimensional symplectic leaf, as is the origin of $\Cone(C,D)$

-- the preimage of $C-D'$ by $\Cone(C,D) \to C$ is an open $2$-dimensional 
symplectic leaf.  

When $C$ is generic, $D = D'$. 
\end{prop}

{\em Proof.} In the proof of Proposition 
\ref{1.1}, $\{\omega,\omega'\}$ has a pole of order $(i+i')\delta_P$, 
so then $\delta_P>0$, this order is $<(i+i'+1)\delta_P$. If we 
view a $k$-differential $\varpi$ as a function on $\Cone(C,D)$, then 
the coefficient of the singularity of order $k\delta_P$ at $P$ should 
be viewed as the value of $\varpi$ at a point of the line of $\Cone(C,D)$
above $P$. So $\{\omega,\omega'\}$ vanishes at $P$ when $\delta_P>0$. 

The elements of $D' - D$ are the points $P$ such that if a section of
$A_i^{(C,D)}$ vanishes at $P$, then it vanishes at $P$ with order $2$. 
Let $(n_1,\ldots,n_g)$ be the Weierstrass sequence of $P$; this sequence 
is defined by the condition that $n_1<\cdots < n_g$, and there exists a 
basis of regular differentials, with zeroes of order $n_1,\ldots,n_g$
at $P$. At a non-Weierstrass point, the sequence is $(0,\ldots,g-1)$. 
Generic curves only have regular Weierstrass points, i.e., with sequence 
$(0,\ldots,g-2,g)$. In both cases, these exist forms $\omega,\omega'\in A_1^{(C)}$, 
such that $\{\omega,\omega'\}$ does not vanish at $P$. 
\hfill \qed \medskip 

If $B$ is an algebra, equipped with a decreasing filtration 
$B = B^{(0)} \supset B^{(1)} \supset \cdots$ (i.e., we have $B^{(i)} B^{(j)} 
\subset B^{(i+j)}$), then its associated graded $\gr(B) = 
\oplus_{i\geq 0} B^{(i)} / B^{(i+1)}$ has a graded ring structure. 
Moreover, if $\gr(B)$ commutative, then it has a natural Poisson structure of 
degree $1$: for $x\in \gr_i(B)$, $y \in \gr_j(B)$, we define $\{x,y\}$ as 
the class of $[\wt x,\wt y]$ in $\gr_{i+j+1}(B)$, where $\wt x,\wt y$ are any 
lifts of $x,y$ in $B^{(i)},B^{(j)}$. We then say that $B$ is a 
quantization of the Poisson algebra $\gr(B)$. 

\medskip 

By a quantization of the Poisson algebra $A^{(D)}$, we therefore understand 
an algebra $B^{(D)}$, together with a decreasing ring filtration, whose
associated graded ring is commutative, and together with an isomorphism 
$\gr(B^{(D)})\to A^{(D)}$ of graded algebras, which is also a Poisson isomorphism. 

\medskip 

The purpose of this paper is to construct a quantization of the 
Poisson algebra (\ref{graded:ring}). Before we explain various forms of 
this construction, let us describe some examples of the Poisson rings 
(\ref{graded:ring}) explicitly in the case $D = 0$ (then the algebra $A^{(D)}$
is simply denoted $A$). We do not know how to quantize the isomorphisms 
$A^{(D_1)} \to A^{(D_2)}$, where $D_1$ and $D_2$ are linearly equivalent. 

\subsection{Explicit form of the Poisson ring $A$ for genus $3,4,5$}

Let us first describe the graded algebra structure of $A$. 
We have $\dimm(A_0) = 1$, and $\dimm(A_1) = g$, where $g$ is the genus of $C$. 
Moreover, the natural map $S^\bullet(A_1) \to A$ is surjective
when $C$ is not hyperelliptic (see \cite{Griffiths:Harris}).
However, the injection $C \hookrightarrow \PP(H^0(C,K)^*)$ is a complete 
intersection only when $g = 3,4,5$ and $C$ is not hyperelliptic. In 
these cases, the kernel of $S^\bullet(A_1) \to A$ is the ideal 
generated by homogeneous elements $Q_1,\ldots,Q_{g-2}$. When 
$g = 3$, $Q_1 = Q$ is homogeneous of degree $4$; 
when $g = 4$, $Q_1,Q_2$ may be taken homogeneous of 
degree $2$ and $3$, and when $g = 5$, $Q_1,Q_2,Q_3$
may all be taken homogeneous of degree $2$ (see \cite{Griffiths:Harris}).

In all these cases, $S^\bullet(A_1)$ may be equipped with a Poisson 
bracket, such that the morphism $S^\bullet(A_1) \to A$ is Poisson; 
in other words, the injection $\Cone(C) \hookrightarrow H^0(C,K)^*$
is a Poisson morphism. The Poisson structure on $S^\bullet(A_1)$
may be described explicitly as follows (see \cite{OR}). 

Let $x_1,\ldots,x_g$ be a basis of $A_1$, then the Poisson structure on 
$S^\bullet(A_1)$ is obtained by the rule 
$$
\{f,g\} = {{df\wedge dg \wedge \chi }\over{\omega_{\on{top}}}}, 
$$  
where $\chi = dQ_1\wedge \cdots \wedge dQ_{g-2}$ and 
$\omega_{\on{top}} = dx_1 \wedge \cdots \wedge dx_g$.
The elements $Q_1,\ldots,Q_{g-2}$ are Poisson central for this 
structure, so there exists a unique Poisson structure on $A$, 
such that $S^\bullet(A_1) \to A$ is Poisson. 
For example, when $g = 3$, the Poisson structure is 
defined by the relations 
$$
\{x_1,x_2\} = \pa_{x_3}Q, \ 
\{x_2,x_3\} = \pa_{x_1}Q, \ 
\{x_3,x_1\} = \pa_{x_2}Q ;  
$$
in general, the brackets have the form 
$\{x_i,x_j\} = P_{ij}(x_1,\ldots,x_g)$, where the 
$P_{ij}$ are homogeneous of degree $3$.

\section{Quantization based on formal pseudodifferential operators} 
\label{sect:pseudo}

\subsection{Outline of the construction}

Our main tool is the general construction of the algebra of formal 
pseudodifferential operators $\PsiDO(R,\pa)$ associated to any differential 
ring $(R,\pa)$. 
We will define the filtered algebra $B$ as an algebra of 
formal pseudodifferential operators on $C$, which are regular on $C$. 
We proceed as follows. To any rational, nonzero vector field $X$ on 
$C$, we associate a filtered algebra $B^{\on{rat}}_X$ of {\it rational}
pseudodifferential operators on $C$. The construction of this algebra 
involves $X$, but we construct canonical isomorphisms
$$
i_{X,Y}^{\on{rat}} : B_X^{\on{rat}} \to B_Y^{\on{rat}}
$$
for any pair $(X,Y)$ of nonzero rational vector fields. One can show that 
$B_X^{\on{rat}}$ is a quantization of the Poisson algebra $A^{\on{rat}}$. 

In Section \ref{rem:indep}, we give a canonical construction of the 
algebra $B^\rat_X$, independent of the choice of a nonzero vector field
$X$. 

If $z$ is a formal variable, $\PsiDO(\CC((z)),{\pa\over{\pa z}})$ is the algebra of formal 
pseudodifferential operators on the formal punctered disc. This algebra contains
the subalgebra  $\PsiDO(\CC[[z]],{\pa\over{\pa z}})$ of operators, regular at the origin.
For any integer $\delta \geq 0$, we also construct an intermediate 
algebra $\PsiDO(\CC[[z]],z^{\delta}{\pa\over{\pa z}})$. 
  
Then for any point $P$ of $C$, let $\cK_P$ be the completed local field of 
$C$ at $P$, and let $\cO_P \subset \cK_P$ be its completed local ring. 
If $z_P$ is a local coordinate at $P$, we have 
$$
\cK_P = \CC((z_P)), \; \cO_P = \CC[[z_P]]. 
$$
We set 
$\PsiDO(\cK_P,z_P) := \PsiDO(\CC((z_P)),{\pa\over{\pa z_P}})$, 
$\PsiDO(\cO_P,z_P) := \PsiDO(\CC[[z_P]],{\pa\over{\pa z_P}})$ and 
$\PsiDO(\cO_P,z_P)^{(\delta_P)} := \PsiDO(\CC[[z_P]],(z_P)^{\delta_P}{\pa\over{\pa z_P}})$.

If $P$ is any point of $C$, Laurent expansion of formal 
pseudodifferential operators at $P$ yields a filtered ring 
morphism $L_P^{z_P} : B_X^{\on{rat}} \to \PsiDO(\cK_P,z_P)_{\leq 0}$. 
Then we define $B_X^{(D)}$ as the preimage of $\prod_{P\in C}
\PsiDO(\cO_P,z_P)^{(\delta_P)}_{\leq 0}$ by the ring morphism 
$$
\prod_{P\in C} L_P^{z_P} : B_X^{\on{rat}} \to 
\prod_{P\in C} \PsiDO(\cK_P,z_P)_{\leq 0} 
$$ 
(the index $\leq 0$ means operators of degree $\leq 0$).
One easily sees that this definition is independent of the 
choice of the collection of local coordintates $(z_P)_{P\in C}$. 
In particular, when $D=0$, $B_X = B_X^{(0)}$ consists of all rational 
pseudodifferential operators on $C$, which are regular at any point of $C$. 
We will prove: 

\begin{thm} \label{thm:pseudo}
1) The canonical isomorphisms $i_{X,Y}^{\on{rat}}$ restrict to 
canonical isomorphisms of filtered algebras 
$$
i_{X,Y} : B^{(D)}_X \to B^{(D)}_Y. 
$$

2) The graded algebra $\gr(B^{(D)}_X)$ is commutative, and as a 
Poisson algebra, it is isomorphic to $A^{(D)}$. 

3) For $D\geq D'$, there are canonical morphisms $B_X^{(D)} \hookrightarrow 
B_X^{(D')}$ of complete filtered algebras, quantizing the inclusion 
$A^{(D)} \hookrightarrow A^{(D')}$. 
\end{thm}

So for each $D$, the algebras $B_X^{(D)}$ are all isomorphic when the vector field
$X$ is changed, and they are quantizations of the Poisson algebra $A^{(D)}$. 

\begin{remark} 
One can prove that if one repeats this construction without restricting it 
to operators of degree $\leq 0$, the resulting algebra is the same as 
$B_X$: all regular pseudodifferential operators on $C$ are of degree $\leq 0$, 
because there are no nonzero sections of $K(D)^{\otimes i}$ for $i<0$ 
(the genus of $C$ is $>1$). 
\end{remark}

\subsection{Details of the construction}

We will first present all details of the construction when 
$D = 0$. So all superscripts $(D)$ will be dropped. 
In Section \ref{outline:D}, we explain the modifications of the 
construction in the case of a general $D$. 

\subsubsection{The algebras $\PsiDO(R,\pa)$}

Let $R$ be a commutative ring with unit and let $\pa$ be a derivation of 
$R$. Following \cite{Adler,Manin}, define $\PsiDO(R,\pa)$ as the space of 
all formal linear combinations $\sum_{i\in\ZZ} a_i D_\pa^i$, 
where for each $i$, $a_i\in R$ and $a_i = 0$ for $i$ large enough. 
$\PsiDO(R,\pa)$ is equipped with the associative product 
$$
(\sum_{i\in\ZZ} a_i D_\pa^i)(\sum_{j\in\ZZ} b_j D_\pa^j)
= \sum_{k\in\ZZ} \big( \sum_{i,j\in\ZZ} \pmatrix i \\ i+j-k \endpmatrix
a_i \pa^{i+j-k}(b_j)\big) D_\pa^k.  
$$
Say that $\sum_{i\in\ZZ} a_i D_\pa^i$ has degree $\leq n$ if $a_i = 0$
when $i>n$, and define $\PsiDO(R,\pa)_{\leq n}$ as the subspace of 
$\PsiDO(R,\pa)$ of all operators of degree $\leq n$. 
Then $\PsiDO(R,\pa)$ is a filtered ring. Its associated graded is 
$R[\xi,\xi^{-1}]$. We will be interested in its subring
$\PsiDO(R,\pa)_{\leq 0}$. It is also filtered, with associated 
graded $R[\xi^{-1}]$. Moreover, both $\PsiDO(R,\pa)$ and $\PsiDO(R,\pa)_{\leq 0}$
are complete for the topology defined
by the family $(\PsiDO(R,\pa)_{\leq -n})_{n = 0,1,2,\ldots}$.  

\subsubsection{Functoriality properties of the rings $\PsiDO(R,\pa)$
and $\PsiDO(R,\pa)_{\leq 0}$}

The following statements are immediate: 

\begin{lemma} \label{lemma:pdo}
1) Let $(R,\pa)$ be a differential ring, and let 
$f\in R^\times$ (so $f$ is an invertible element of $R$). 
Set $\pa' = f\pa$, then $\pa'$ is a derivation of $R$. 
We have for any $i$, $(f^{-1}D_{\pa'})^i\in\PsiDO(R,\pa')_{\leq i}$, 
so if $(a_i)_{i\in\ZZ}$ is a sequence of elements of $R$, 
such that $a_i = 0$ for $i$ large enough, the sequence 
$\sum_{i\in \ZZ} a_i (f^{-1} D_{\pa'})^i$ converges in $\PsiDO(R,\pa')$. 
Then there is a unique isomorphism
$$
i_{\pa,\pa'} : \PsiDO(R,\pa) \to \PsiDO(R,\pa')
$$
of complete filtered algebras, taking each series 
$\sum_{i\in\ZZ} a_i D_\pa^{i}$ to $\sum_{i\in \ZZ} a_i (f^{-1} D_{\pa'})^i$. 
We have then $i_{\pa',\pa''} \circ i_{\pa,\pa'} = i_{\pa,\pa''}$. 

2) Let $\mu : (R,\pa_R) \to (S,\pa_S)$ be a morphism of differential 
rings (i.e., $\mu$ is a ring morphism and $\pa_S \circ \mu = \mu \circ \pa_R$). 
Then there is a unique morphism 
$$
\PsiDO(\mu) : \PsiDO(R,\pa_R) \to \PsiDO(S,\pa_S), 
$$
taking each $\sum_{i\in\ZZ} a_i D_{\pa_R}^i$ to $\sum_{i\in\ZZ}
\mu(a_i) D_{\pa_S}^i$. $\PsiDO(\mu)$ is a morphism of 
complete filtered algebras and we have $\PsiDO(\nu\circ \mu) = 
\PsiDO(\nu) \circ \PsiDO(\mu)$ for any morphism $\nu : (S,\pa_S)
\to (T,\pa_T)$ of differential rings. In other words, 
$$
(R,\pa) \mapsto \PsiDO(R,\pa)
$$  
is a functor from the category of differential rings to that of
filtered complete algebras. 
\end{lemma}

\subsubsection{Construction of $B_X$}

Let $C$ be a curve, and let $\CC(C)$ be its field of rational functions. 
Let $X$ be a nonzero rational vector field on $C$; $X$ may be viewed as 
a nonzero derivation of $\CC(C)$. We set 
$$
B_X^{\on{rat}} = \PsiDO(\CC(C),X)_{\leq 0}. 
$$
If $Y$ is another nonzero vector field on $C$, then there exists 
a unique $f\in\CC(C)^\times$, such that $Y = fX$. Applying 
Lemma \ref{lemma:pdo}, 1), we get an isomorphism 
$$
i_{X,Y}^{\on{rat}} : B_X^{\on{rat}} \to B_Y^{\on{rat}}
$$
of complete filtered rings. 

On the other hand, if $P\in C$, then for any local coordinate $z_P$ at $P$, 
${\pa\over{\pa z_P}}$ is a derivation of $\cK_P$, preserving $\cO_P$. 
we set 
$$
\PsiDO(\cK_P,z_P) := \PsiDO(\cK_P,{\pa\over{\pa z_P}})
\; \on{and} \; 
\PsiDO(\cO_P,z_P) := \PsiDO(\cO_P,{\pa\over{\pa z_P}}). 
$$
By functoriality, we have then an inclusion 
$\PsiDO(\cO_P,z_P) \subset \PsiDO(\cK_P,z_P)$. Moreover, if 
$z'_P$ is another local coordinate at $P$, the derivations 
${\pa\over{\pa z_P}}$ and ${\pa\over{\pa z'_P}}$ are related by 
${\pa\over{\pa z_P}} = \varphi {\pa\over{\pa z'_P}}$, 
where $\varphi$ belongs to 
$\cO_P^\times$, so Lemma \ref{lemma:pdo}, 1), says that there is an 
isomophism $i_{z_P,z'_P} : \PsiDO(\cK_P,z_P) \to \PsiDO(\cK_p,z'_P)$
of complete filtered algebras, restricting to an isomorphism 
$\PsiDO(\cO_P,z_P) \to \PsiDO(\cO_p,z'_P)$, and such that 
$i_{z'_P,z''_P} \circ i_{z_P,z'_P}  = i_{z_P,z''_P}$. 

Let us now define the Laurent expansion morphism 
$$
L_P^{z_P} : B_X^{\on{rat}} \to \PsiDO(\cK_P,z_P)_{\leq 0}. 
$$
Since $X$ is a nonzero vector field, its local expansion at $P$
is $X = X(z_P) {\pa\over{\pa z_P}}$, with $X(z_P)\in\cK_P^\times$.
The Laurent expansion map 
$$
\ell_P : \CC(C) \to \cK_P
$$
therefore induces a differential ring morphism $(\CC(C),X) \to 
(\cK_P,X(z_P){\pa\over{\pa z_P}})$, and so a morphism 
$$
\PsiDO(\ell_P) : B_X^{\on{rat}} \to \PsiDO(\cK_P,X(z_P){\pa\over{\pa z_P}})_{\leq 0}. 
$$
Composing it with the isomorphism 
$$
i_{X(z_P){\pa\over{\pa z_P}},{\pa\over{\pa z_P}}} : 
\PsiDO(\cK_P,X(z_P){\pa\over{\pa z_P}})_{\leq 0}
\to \PsiDO(\cK_P,{\pa\over{\pa z_P}})_{\leq 0}
= \PsiDO(\cK_P,z_P)_{\leq 0} , 
$$
we get a filtered ring morphism 
$$
L_P^{z_P} : B_X^{\on{rat}} \to \PsiDO(\cK_P,z_P)_{\leq 0}. 
$$
Finally, let us prove that the preimage by 
$$
\prod_{P\in C} L_P^{z_P} : B_X^{\on{rat}} \to \prod_{P\in C}
\PsiDO(\cK_P,z_P)_{\leq 0}
$$
of $\prod_{P\in C}\PsiDO(\cO_P,z_P)_{\leq 0}$ is independent of the 
choice of the local coordinates $(z_P)_{P\in C}$:  if $(z'_P)_{P\in C}$
is any other choice of local coordinates, then 
\begin{equation} \label{vfields:Laurent}
L_P^{z'_P} = i_{{\pa\over{\pa z_P}},{\pa\over{\pa z'_P}}} \circ 
L_P^{z_P},
\end{equation}
and 
$$
\PsiDO(\cO_P,z'_P)_{\leq 0} =  i_{{\pa\over{\pa z_P}},{\pa\over{\pa z'_P}}} 
(\PsiDO(\cO_P,z'_P)_{\leq 0} ), 
$$
so
\begin{align*}
\big(\prod_{P\in C} L_P^{z'_P}\big)^{-1}
\big(\prod_{P\in C}\PsiDO(\cO_P,z'_P)_{\leq 0}\big) 
& = \big(\prod_{P\in C} i_{{\pa\over{\pa z_P}},{\pa\over{\pa z'_P}}}
\circ L_P^{z_P}\big)^{-1} \big(\prod_{P\in C}\PsiDO(\cO_P,z'_P)_{\leq 0}\big) 
\\ & = \big(\prod_{P\in C} L_P^{z_P}\big)^{-1}
\big(\prod_{P\in C}\PsiDO(\cO_P,z_P)_{\leq 0}\big) . 
\end{align*}

\subsubsection{Vector field-independent construction of the algebras 
$B_X^\rat$} \label{rem:indep}

Let us define $\on{DO}(\CC(C))$ as the algebra of all 
rational differential operators on $C$. So $\on{DO}(\CC(C))$
is generated by $i(f)\in \CC(C)$, $D_X$, where $X\in \Der(\CC(C))$, 
and relations 
$$
i(fg) = i(f)i(g), \quad i(\la) = \la, \quad [D_X,i(f)] = i(X(f)), 
$$
\begin{equation} \label{def:eqs}
[D_X,D_Y] = D_{[X,Y]}, \quad D_{fX} = i(f)D_X, 
\end{equation}
for $f,g\in\CC(C)$, $X,Y\in \Der(\CC(C))$, $\la\in\CC$. We will
denote $i(f)$ simply by $f$, for $f\in\CC(C)$. 

$\on{DO}(\CC(C))$ may be localized with respect to the family of 
all $D_X$, where $X$ are all nonzero rational vector fields. 
The last of relations (\ref{def:eqs}), together with the fact that 
$\Der(\CC(C))$ is a 1-dimensional $\CC(C)$-vector space, implies that 
the localization of $\on{DO}(\CC(C))$ w.r.t.\  any $D_X$, 
$X\in\Der(\CC(C)) - \{0\}$, coincides with its localization w.r.t.\   
the family of all such $D_X$. We denote by $B^\rat$ the completion of this 
localized algebra w.r.t.\ the degree of formal pseudodifferential operators. 

Then $B^\rat$ contains $\on{DO}(\CC(C))$ as a subalgebra, as well as the 
additional generators $(D_X)^{-1}$, $X\in \Der(\CC(C)) - \{0\}$. They
satisfy, in particular, the relations 
$$
(D_{fX})^{-1} = (D_X)^{-1} f^{-1}, 
$$
for $f\in \CC(C)^\times$ and $X\in \Der(\CC(C)) - \{0\}$. If $X$
is any nonzero vector field, the natural map
$$
i_X : B_X^\rat \to B^\rat 
$$
is therefore an isomorphism. The map $i^\rat_{X,Y} : B_X^\rat \to B_Y^\rat$ 
then coincides with $(i_Y)^{-1} \circ i_X$. 

\subsubsection{Proof of Theorem \ref{thm:pseudo}} \label{sect:proof}

Let us prove the first part of Theorem \ref{thm:pseudo}.
Let us emphasize the dependence of $L_P^{z_P}$ in $X$
by denoting it 
$$
L_P^{X,z_P} : B_X^{\on{rat}} \to \PsiDO(\cK_P,z_P)_{\leq 0}. 
$$
Then we have 
$$
B_X = \big( \prod_{P\in C} L_P^{X,z_P}\big)^{-1} 
(\prod_{P\in C} \PsiDO(\cO_P,z_P)). 
$$ 
Now the composed map 
$$
B_X^{\on{rat}} \stackrel{i_{X,Y}}{\to} B_Y^{\on{rat}}
\stackrel{L_P^{Y,z_P}}{\to} \PsiDO(\cK_P,z_P)_{\leq 0}
$$
coincides with $L_P^{Y,z_P}$. So 
\begin{align*}
B_X & =
\big( \prod_{P\in C} L_P^{Y,z_P} \circ i^{\on{rat}}_{X,Y}\big)^{-1} \big( \prod_{P\in C}
\PsiDO(\cO_P,z_P)\big) 
\\ & = (i^{\on{rat}}_{X,Y})^{-1} \Big( \big( \prod_{P\in C} L_P^{Y,z_P} \big)^{-1} 
\big( \prod_{P\in C} \PsiDO(\cO_P,z_P)\big) \Big) 
\\ & = (i^{\on{rat}}_{X,Y})^{-1}(B_Y), 
\end{align*}
so $B_X = (i^{\on{rat}}_{X,Y})^{-1}(B_Y)$. Since $i^{\on{rat}}_{X,Y} : B_X^{\on{rat}} 
\to B_Y^{\on{rat}}$ is an isomorphism of complete filtered algebras, 
it restricts to an isomorphism $i_{X,Y} : B_X \to B_Y$ of complete filtered algebras. 

Let us now prove the second part of Theorem \ref{thm:pseudo}.
We will define a filtration on $B_X$; then we will construct 
a graded linear map 
$$
\la_{\on{reg}} : \gr(B_X) \to A; 
$$
we will prove that if the genus of $C$ is $>1$, $\la_{\on{reg}}$ is a linear 
isomorphism, and finally that it is an isomorphism of Poisson algebras. 

\medskip 
{\it (a) Filtration on $B_X$.} 
We set $(B_X^{\on{rat}})^i = \PsiDO(\CC(C),X)_{\leq -i}$, 
and 
$$
(B_X)^i = B_X \cap (B_X^{\on{rat}})^i. 
$$
So $(B_X)^i$ consists of all regular pseudodifferential 
operators on $C$ of order $\leq -i$. 

\medskip 
{\it (b) The map $(B_X)^i/(B_X)^{i+1} \to A_i$.} 
The natural map $(B_X)^i/(B_X)^{i+1} \to (B^{\on{rat}}_X)^i/(B^{\on{rat}}_X)^{i+1}$
is injective, because $(B_X)^i \cap (B_X^{\on{rat}})^{i+1} = 
B_X \cap (B_X^{\on{rat}})^{i+1} = (B_X)^{i+1}$. Moreover, there is a
linear isomorphism 
\begin{equation} \label{la:i:rat} 
\la^{(i)}_{\on{rat}} : (B_X^{\on{rat}})^i / (B_X^{\on{rat}})^{i+1} 
\to \{\on{rational}\  i\on{-differentials\  on\ } C\},  
\end{equation} 
taking the class of $\sum_{j\geq i} a_j (D_X)^{-j}$ to $a_i \al^i$, 
where $\al$ is the rational differential inverse to $X$. 
We will prove 

\begin{lemma}
The restriction of $\la^{(i)}_{\on{reg}}$ maps 
$(B_X)^i/(B_X)^{i+1}$ to $A_i = H^0(C,K^{\otimes i}) \subset  
\{$rational $i$-differentials on $C\}$. 
\end{lemma}

{\em Proof of Lemma.} For any $P\in C$, $L_P^{z_P}$ induces a linear 
map 
$$
(B_X^{\on{rat}})^i / (B_X^{\on{rat}})^{i+1} \to \PsiDO(\cK_P,z_P)_{\leq -i} / 
\PsiDO(\cK_P,z_P)_{\leq -i-1} ; 
$$
it restricts to a linear map 
$$
(B_X)^i / (B_X)^{i+1} \to \PsiDO(\cO_P,z_P)_{\leq -i} / 
\PsiDO(\cO_P,z_P)_{\leq -i-1} .  
$$
Now we have a linear isomorphism
$$
\la_P^{(i)} : \PsiDO(\cK_P,z_P)_{\leq -i} / \PsiDO(\cK_P,z_P)_{\leq -i-1}
\to \CC((z_P)) (dz_P)^{\otimes i}
$$
restricting to an isomorphism  
$$
\PsiDO(\cO_P,z_P)_{\leq -i} / \PsiDO(\cO_P,z_P)_{\leq -i-1}
\to \CC[[z_P]] (dz_P)^{\otimes i}
$$
and taking the class of $\sum_{j\geq i} b_j \pa^{-j}$ to the class
of $\wt b_i (dz_P)^{\otimes i}$, where $\wt b_i$ corresponds to 
$b_i$ under $\cK_p = \CC((z_P))$, and the diagram 
$$
\begin{array}{ccc}
(B_X^{\on{rat}})^i / (B_X^{\on{rat}})^{i+1} & @>\la^{(i)}_{\on{rat}}>> & 
\{\on{rational\ }i\on{-differentials}\}
\\ 
\downarrow{L_P^{z_P}} & \; & \downarrow{\ell_P^{(i)}}
\\ 
\PsiDO(\cK_P,z_P)_{\leq -i} / \PsiDO(\cK_P,z_P)_{\leq -i-1}
 & @>\la_P^{(i)}>> & \CC((z_P)) (dz_P)^{\otimes i}
\end{array}
$$
is commutative (the right vertical arrow $\ell_P^{(i)}$ is the Laurent expansion 
on $i$-differentials at $P$). Then $L_P^{z_P}$ maps  
$(B_X)^i / (B_X)^{i+1}$ to  
$$
\PsiDO(\cO_P,z_P)_{\leq -i} / \PsiDO(\cO_P,z_P)_{\leq -i-1},
$$
so $\la_P^{(i)} \circ L_P^{(z_P)}$
maps $(B_X)^i/(B_X)^{i+1}$ to $\CC[[z_P]](dz_P)^{\otimes i}$. Therefore 
$\la^{(i)}_{\on{reg}}((B_X)^i/(B_X)^{i+1})$ is contained in the space of 
rational differentials on $C$, which are regular at each point of $C$; 
this space is precisely $H^0(C,K^{\otimes i}) = A_i$. 
\hfill \qed \medskip 

Being the restriction of an injective map, the map 
$\la^{(i)}_{\on{reg}} : (B_X)^i / (B_X)^{i+1} \to A_i$ induced 
by $\la^{(i)}_{\on{rat}}$ is injective. We now prove: 

{\it (c) The map $\la^{(i)}_{\on{reg}} : (B_X)^i / (B_X)^{i+1} \to A_i$
is surjective.} 

Let $(a_j)_{j = i,i+1,\ldots}$ be a collection of elements of 
$\CC(C)$; let us write the necessary and sufficient conditions 
for $\sum_{j\geq i} a_j D_X^{-j}$ to be a regular pseudodifferential 
operator. For simplicity, we will assume that the form $\al = X^{-1}$
has no pole and $2g-2$ distinct zeroes $Q_1,\ldots,Q_{2g-2}$; so the 
vector field $X$ is nowhere vanishing and has simple poles at $Q_1,\ldots,
Q_{2g-2}$. Let $z_\al$ be a local coordinate at $Q_\al$. Then we have a 
local expansion at $Q_\al$
$$
X = \big( {{c_\al}\over{z_\al}} + \on{element\ of\ }\CC[[z_\al]]\big) 
{\pa\over{\pa z_\al}},  
$$
where $c_\al\in\CC^\times$; so we have local expansions 
\begin{align*}
(D_X)^{-j} = & \la_{j,j}^{(\al)} (z_\al)^j (D_{\pa/\pa z_\al})^{-j}
+ \la_{j,j+1}^{(\al)} (z_\al)^{j-1} (D_{\pa/\pa z_\al})^{-j-1}
+ \cdots + \la_{j,2j-1}^{(\al)} z_\al (D_{\pa/\pa z_\al})^{-2j+1}
\\ & 
+\sum_{\ell\geq 0} \la_{j,2j+\ell}^{(\al)} (D_{\pa/\pa z_\al})^{-2j-\ell} ,  
\end{align*}
where $\la_{j,j}^{(\al)},\ldots,\la_{j,2j-1}^{(\al)}\in \CC[[z_\al]]^\times$, 
and $\la_{j,2j}^{(\al)}, \la_{j,2j+1}^{(\al)},\ldots \in \CC[[z_\al]]$
(the constant terms of $\la_{j,j}^{(\al)},\ldots,\la_{j,2j-1}^{(\al)}$ may
be computed explicitly using binomial coefficients). Recall that $D_X$
is the generator of $B_X^{\on{rat}}$ corresponding to the vector field
$X$, and $D_{\pa/\la z_\al}$ is the generator of $\PsiDO(\cK_{Q_\al},z_\al)$
corresponding to the vector field ${\pa\over{\pa z_\al}}$. 

So the necessary and sufficient conditions on $(a_j)_{j = i,i+1,\ldots}$
are: 

(a) $a_j\in \CC(C)$, and each $a_j$ is regular outside
$\{P_1,\ldots,P_{2g-2}\}$; 

(b) (local conditions are each $Q_\al$, $\al = 1,\ldots,2g-2$) 
let us denote by $a_j^{(\al)}$ the element of $\CC((z_\al))$, 
obtained as the Laurent expansion of $a_j$ at $Q_\al$, then the 
formal series 
$$
\la_{i,i}^{(\al)} (z_\al)^i a_i^{(\al)}, \; \;
\la_{i+1,i+1}^{(\al)} (z_\al)^{i+1} a_{i+1}^{(\al)}
+ \la_{i,i+1}^{(\al)} (z_\al)^{i-1} a_{i}^{(\al)},
$$
$$
\la_{i+2,i+2}^{(\al)} (z_\al)^{i+2} a_{i+2}^{(\al)}
+ \la_{i+1,i+1}^{(\al)} (z_\al)^{i} a_{i+1}^{(\al)}
+ \la_{i,i+2}^{(\al)} (z_\al)^{i-2} a_{i}^{(\al)},
\; \on{etc.}
$$
should all be regular. 

This means that the formal series $(a_j^{(\al)})_{j = i,i+1,\ldots}$
should have the expansions: 
$$
a_i^{(\al)} = \al_{i,i} (z_\al)^{-i} + \al_{i,i-1} (z_\al)^{-i+1} + \cdots, 
$$
$$
a_{i+1}^{(\al)} = A_{i+1,i+2}(a_i^{(\al)}) (z_\al)^{-i-2}
+ \al_{i+1,i+1} (z_\al)^{-i-1} + \al_{i+1,i} (z_\al)^{-i} + \cdots, 
$$
\begin{align*}
a_{i+2}^{(\al)} = & A_{i+2,i+4}(a_{i+1}^{(\al)}) (z_\al)^{-i-4}
+ A_{i+2,i+3}(a_{i+1}^{(\al)}) (z_\al)^{-i-3}
\\ & 
+ \al_{i+2,i+2} (z_\al)^{-i-2} + \al_{i+2,i+1} (z_\al)^{-i-1} + \cdots, 
\end{align*}
where the $\al_{k,l}$ are arbitrary complex numbers, and the 
$f\mapsto A_{k,l}(f)$ are certain linear forms on $(z_\al)^{-k}\CC[[z_\al]]$. 

These conditions can be translated as follows: 

(1) $a_i X^{-i}\in H^0(C,K^{\otimes i})$; 

(2) $a_{i+1}$ belongs to a (possibly empty) affine space over 
$H^0(C,K^{\otimes i+2})$, depending on $a_i$; 

(3) $a_{i+2}$ belongs to a (possibly empty) affine space over 
$H^0(C,K^{\otimes i+3})$, depending on $a_{i+1}$, etc. 

We now prove that these affine spaces are all nonempty, and we 
describe the set of all possible $(a_j)_{j\geq i}$. 

\begin{lemma} \label{lemma:2.3}
Let $i\geq 2$ and $j\geq 0$. Define $D_{\on{can}}$ as the 
divisor $Q_1 + \cdots + Q_{2g-2}$. Identify $H^0(C,K^{\otimes i})$
with the space $\{f\in \CC(C)\¦ (f)\geq -i D_{\on{can}}\}$, and 
$H^0(C,K^{\otimes i}(jD_{\on{can}}))$ with $\{f\in \CC(C)\¦ (f)\geq -(i+j) D_{\on{can}}\}$. 
Then we have $H^0(C,K^{\otimes i}) \subset H^0(C,K^{\otimes i}(jD_{\on{can}}))$. 
Moreover, for $\al = 1,\ldots,2g-2$, and $k = 1,\ldots,j$, define
linear forms 
$$
\phi_{\al,k} : H^0(C,K^{\otimes i}(jD_{\on{can}})) \to \CC
$$
by the condition that the local expansion of $f$ at $Q_\al$ is 
$$
f\in \sum_{k = 1}^j \phi_{\al,k}(f) (z_\al)^{-i-k} + (z_\al)^{-i}
\CC[[z_\al]]. 
$$ 
Then the sequence 
\begin{equation} \label{exact}
0 \to H^0(C,K^{\otimes i}) \to H^0(C,K^{\otimes i}(jD_{\on{can}})) 
\stackrel{\oplus_{\al,k} \phi_{\al,k}}{\to} \CC^{(2g-2)j} \to 0
\end{equation}
is exact. 
\end{lemma}

{\em Proof of Lemma.} This follows from the fact that if $(R_1,\ldots,R_\ell)$
is any sequence of points of $C$, then 
$$
h^0(C,K^{\otimes i}(R_1 + \cdots + R_{\ell}))
= h^0(C,K^{\otimes i}(R_1 + \cdots + R_{\ell-1})) +1, 
$$ 
because the degree of $K^{\otimes i}(R_1 + \cdots + R_{\ell -1})$
is $\geq 2(g-1)$. 
\hfill \qed \medskip 

For any pair $(i,j)$, let us choose a section $\sigma_{i,j}$
of the exact sequence (\ref{exact}). So $\sigma_{i,j}$ is a linear 
map 
$$
\sigma_{i,j} : H^0(C,K^{\otimes i}) \to H^0(C,K^{\otimes i}(jD_{\on{can}})), 
$$ 
such that if $f = \sigma_{i,j}((\la_{\al,k})_{\al,k})$, then 
for each $(\al,k)$, we have $\phi_{\al,k}(f) = \la_{\al,k}$. 

For any $\omega\in H^0(C,K^{\otimes i})$, we set 
$$
\sigma(\omega) = (a_i(\omega),a_{i+1}(\omega),\ldots),  
$$
where 
$$
a_i(\omega) = \omega X^{-i}, 
$$
$$
a_{i+1}(\omega) = \sigma_{i+1,1}( \big(  A_{i+1,i+2}(a_i(\omega)^{(\al}))
 \big)_{\al = 1,\ldots,2g-2})), 
$$
$$
a_{i+2}(\omega) = \sigma_{i+2,2}( \big( A_{i+2,i+4}(a_{i+1}(\omega)^{(\al)}), 
A_{i+3,i+4}(a_{i+1}(\omega)^{(\al)}) \big)_{\al = 1,\ldots,2g-2})), 
$$
etc. Then $\sigma$ is a linear map 
$$
\sigma : H^0(C,K^{\otimes i}) \to B_X^i;  
$$
it is a section of the canonical projection $\la^{(i)}_{\on{reg}} : 
B_X^{(i)} \to H^0(C,K^{\otimes i})$. This proves that $\la^{(i)}_{\on{reg}}$ 
is surjective. 

{\it (d) The map $\la : \gr(B_X) \to A$ is an isomorphism of Poisson 
algebras.}

There is a unique Poisson structure on $\CC(C)[\xi^{-1}]$, such that 
$\{f,g\} = 0$ and 
$$
\{\xi^{-1},f\} = - X(f) \xi^{-2} 
$$
for $f,g\in\CC(C)$. Then the map $\la_{\on{rat}} : 
\gr(B_X^{\on{rat}}) \to \CC(C)[\xi^{-1}]$ is an 
isomorphism of Poisson algebras. Moreover, there is a unique 
inclusion  $A \hookrightarrow \CC(C)[\xi^{-1}]$, taking $\omega\in 
H^0(C,K^{\otimes i})$ to $(\omega \al^{-i}) \xi^{-i}$ (recall that 
$\omega\al^{-i}$ belongs to $\CC(C)$). This inclusion is a morphism
of Poisson algebras. Then we have a commuting diagram 
$$
\begin{array}{ccc}
\gr(B_X^{\on{rat}}) & @>\la_{\on{rat}}>> & \CC(C)[\xi^{-1}]
\\ 
\uparrow & \; & \uparrow
\\ 
\gr(B_X) & @>\la_{\on{reg}}>> & A
\end{array}
$$
Since all maps in this diagram except possibly $\la_{\on{reg}}$
are Poisson algebra morphisms, and since the vertical arrows are
injective, $\la_{\on{reg}}$ is also a Poisson morphism. This ends
the proof of Theorem  \ref{thm:pseudo}, in the case $D =0$. \hfill \qed \medskip  

\subsubsection{The case of nonzero divisor $D$} \label{outline:D}

We already defined the algebra $B_X^{(D)}$, using the vector 
field $X$ and the collection of local coordinates 
$(z_P)_{P\in C}$. We first prove:  

\begin{lemma}
$B_X^{(D)}$ is independent of the choice of $(z_P)_{P\in C}$. 
\end{lemma}

{\em Proof.}
If $z_P$ and $z'_P$ are local coordinates at $P\in C$, we have 
an isomorphism $i_{z_P,z'_P} : \PsiDO(\cK_P,z_P) \to \PsiDO(\cK_P,z'_P)$. 
After composing it with the isomorphisms 
$i_{ {\pa\over{\pa z_P}} , (z_P)^{\delta_P} {\pa\over{\pa z_P}} } : 
\PsiDO(\cK_P,z_P) \to \PsiDO(\cK_P,(z_P)^{\delta_P}{\pa\over{\pa z_P}})$
and  the inverse of 
$i_{ {\pa\over{\pa z'_P}} , (z'_P)^{\delta_P} {\pa\over{\pa z'_P}} } : 
\PsiDO(\cK_P,z'_P) \to \PsiDO(\cK_P,(z'_P)^{\delta_P}{\pa\over{\pa z'_P}})$, 
we get the isomorphism 
$$
i_{ (z_P)^{\delta_P} {\pa\over{\pa z_P}} , (z'_P)^{\delta_P} {\pa\over{\pa z'_P}} } : 
\PsiDO(\cK_P,(z_P)^{\delta_P}{\pa\over{\pa z_P}})
\to \PsiDO(\cK_P,(z'_P)^{\delta_P}{\pa\over{\pa z'_P}}). 
$$
Now there exists $\varphi\in\cO_P^\times$, such that 
$(z_P)^{\delta_P} {\pa\over{\pa z_P}} = \varphi \cdot 
(z'_P)^{\delta_P} {\pa\over{\pa z'_P}}$, 
so $i_{ (z_P)^{\delta_P} {\pa\over{\pa z_P}} , (z'_P)^{\delta_P} {\pa\over{\pa z'_P}} }$
restricts to an isomorphism 
$$
\PsiDO(\cO_P,z_P)^{(\delta_P)} \to \PsiDO(\cO_P,z'_P)^{(\delta_P)}. 
$$
One uses this isomorphism in the same way as above to show that 
the algebra $B_X^{(D)}$ is independent on the choice of 
$(z_P)_{P\in C}$. \hfill \qed \medskip 

The behavior of $B_X^{(D)}$ with respect to changes of the vector field $X$
is proved as above. 

The filtration of $B_X^{(D)}$ is defined by 
$(B_X^{(D)})^i = B_X^{(D)} \cap (B_X^{\on{rat}})^i$. Then we prove:  

\begin{lemma}
The restriction of the map $\la^{(i)}_{\on{rat}}$
defined by (\ref{la:i:rat}) to  $(B_X^{(D)})^i  / (B_X^{(D)})^{i+1}$
maps to $(A^{(D)})_i \subset A^{\on{rat}}_i$. So $\la_{\on{rat}}$
induces a Poisson morphism $\gr(B_X^{(D)}) \to A^{(D)}$. 
\end{lemma}

{\em Proof.} Any element of $\PsiDO(\cO_P,
(z_P)^{\delta_P}{\pa\over{\pa z_P}})_{\leq -i}$
has the expansion 
$$
\sum_{j\geq i} a_j (D_{\pa/\pa z_P})^{-j},
$$
with $a_j\in (z_P)^{-j\delta_P}\CC[[z_P]]$ for any $j\geq i$. This 
implies that $a_i \al^i\in H^0(C,K(D)^{\otimes i})$. 
\hfill \qed \medskip 

The other statements are proved as above, in particular, the analogue of 
Lemma \ref{lemma:2.3} holds because $\on{deg}(K(D)) \geq \on{deg}(K)$. 
\hfill \qed \medskip 

\subsection{Twisting by generalized line bundles} \label{line}

If $\ell\in\CC$, there is a unique continuous automorphism of 
$\PsiDO(\CC((z)),{\pa\over{\pa z}})$, taking $D_{\pa / \pa z}$
to $D_{\pa / \pa z} - {\la\over z}$ and leaving $z$ fixed. 
We denote it by $T \mapsto z^\ell T z^{-\ell}$. 

We denote by $\CC C$ the group of all formal linear combinations 
$\sum_{P\in C} \la_P P$, where all $\la_P$ but a finite number 
are zero. We have a natural group morphism $\ZZ C \to \CC C$. 
Moreover, the residue map is a group morphism $\res : 
\CC(C)^\times \to \ZZ C$. The Picard group of $C$ is defined as 
$\Pic(C) = \ZZ C / \res(\CC(C)^\times)$. Then there is an injection 
$\Pic(C) \hookrightarrow \CC C / \res(\CC(C)^\times)$ induced 
by $\ZZ C \hookrightarrow \CC C$. We call elements of $\CC C$
"generalized divisors" and elements of $\CC C / \res(\CC(C)^\times)$ 
"generalized line bundles".  

Let $\la = \sum_{P\in C} \la_P P$ be a generalized divisor. 
One can define an algebra $B_{X}^{(C,D),\la}$ of twisted  pseudodifferential 
operators as follows
$$
B_X^{(C,D),\la} = \big( \prod_{P\in C} L_P^{z_P}\big)^{-1}
\big( \prod_{P\in C}
(z_P)^{\la_P} \PsiDO(\CC[[z_P]], (z_P)^{\delta_P}{\pa\over{\pa z_P}})
(z_P)^{-\la_P}  \big) . 
$$
Conjugation by a rational function sets up an isomorphism between 
$B_X^{(C,D),\la}$ and $B_X^{(C,D),\la'}$, for $\la,\la'$ linearly equivalent 
generalized divisors (i.e., differing by an element of $\res(\CC(C)^\times)$). 
On the other hand, one can repeat the proof of Theorem \ref{thm:pseudo} to 
prove that the $B_X^{(C,D),\la}$ are quantizations of $A^{(C,D)}$ 
for any $\la$.


\subsection{Functoriality in $(C,D)$} \label{sect:functoriality}

In this section, we emphazise the dependence of the algebras 
$A,A^{(D)},B_X,B_X^{(D)}$ in the curve $C$ by denoting them 
$A^{(C)}$, $A^{(C,D)}$, $B^{(C)}_X$, $B_X^{(C,D)}$. 

\subsubsection{The Poisson algebras}

Let $\varphi : C \to C'$ be a (possibily ramified) covering. So 
$\varphi$ gives rise to an 
inclusion of fields $\varphi^* : \CC(C') \hookrightarrow \CC(C)$. 
Then $\vp$ induces morphisms $\vp^* : H^0(C',(K')^{\otimes i}) \to 
H^0(C,K^{\otimes i})$ (here $K'$ is the canonical bundle of $C'$), 
and therefore an algebra morphism
\begin{equation} \label{phi:no:div}
\vp^*_{\class} : A^{(C')} \to A^{(C)}. 
\end{equation}
The maps $\vp^*$ extend to maps between spaces of 
rational $i$-differentials. For any $f'\in \CC(C')$, we have in 
particular $d(\vp^*(f')) = \vp^*(df')$. It follows that for any rational 
differential $\al'$ on $C'$, we have $\vp^*(\nabla^{\al'}(\omega')) = 
\nabla^{\vp^*(\al')}(\vp^*(\omega'))$. It follows that $\vp^*_\class$
is a morphism of Poisson algebras. 

Let $D' = \sum_{P'\in C'} \delta'_{P'} P'$ be an effective divisor on 
$C'$, and set 
$$
\vp^{-1}(D') = \sum_{P\in C} \big( \delta_{\vp(P)} \nu_P + (1-\nu_P)\big) P, 
$$
where $\nu_P$ is the ramification index of $f$ at $P\in C$
(it is $1$ for all but finitely many $P$). Set $D = \vp^{-1}(D')$, then 
$D$ is an effective divisor of $C$. Then $\vp$ induces a morphism 
\begin{equation} \label{phi:div}
\vp^*_{\class} : A^{(C',D')} \to A^{(C,D)} 
\end{equation}
of graded algebras and of Poisson algebras. 

\subsubsection{Quantization of the morphisms $\vp^*_\class$}

Let $X'$ be a rational, nonzero vector field on $C'$, 
let $\al' = (X')^{-1}$ be the rational differential on $C'$ inverse to $X'$; 
let us set $\al = \vp^*(\al')$ and $X = \al^{-1}$. So $X$ is a 
rational, nonzero vector field on $C$. 
We will now show: 
\begin{prop} \label{prop:funct}
There exists a morphism 
$$
\vp^*_\rat : B_{X'}^{C',\rat} \to B_{X}^{C,\rat}
$$ 
of complete filtered algebras. It induces morphisms
$$
\vp^*_\pseudo :  B_{X'}^{(C')} \to B_{X}^{(C)}
$$
and 
$$
\vp^*_\pseudo :  B_{X'}^{(C',D')} \to B_{X}^{(C,D)}
$$
of complete filtered algebras, quantizing the morphisms 
$\vp^*_\class$. 
\end{prop}

{\em Proof.}

\medskip \noindent

{\it (a) Construction of $\vp^*_\rat : B_{X'}^{C',\rat} \to 
B_{X}^{C,\rat}$.} The map 
$\vp^* : (\CC(C'),X') \to (\CC(C),X)$ is a morphism of differential 
rings. Indeed, we have, for $f'\in \CC(C')$, 
$$
X(\vp^*(f')) = {{d(\vp^*(f'))}\over{\al}} = {{d(\vp^*(f'))}\over{\vp^*(\al')}}
= \vp^*({{df'}\over{\al'}}) = \vp^*(X'(f')).  
$$
So $\vp^*$ induces an algebra map 
$$
\PsiDO(\vp^*)_{\leq 0} : \PsiDO(\CC(C'),X')_{\leq 0}  \to \PsiDO(\CC(C),X)_{\leq 0} , 
$$
that is an algebra map $\PsiDO(\vp^*)_{\leq 0} : B_{X'}^{C',\rat} \to B_X^{C,\rat}$. 

\medskip \noindent

{\it (b) $\vp^*_\rat(B_{X'}^{(C')}) \subset B_X^{(C)}$.}
Let $P\in C$, and let us set $P' = \vp^*(P)$. Let $\nu_P$
be the ramification index of $\vp$ at $P$. Then if $z_P$, 
$z'_{P'}$ are local coordinates at $P,P'$, we have 
$\vp^*(z'_{P'}) = \la \cdot (z_P)^{\nu_P}$, where 
$\la\in \CC[[z_P]]^\times$. Then there is a natural morphism 
$$
\PsiDO(\cK_{P'},z'_{P'})_{\leq 0} \to \PsiDO(\cK_{P},z_{P})_{\leq 0},  
$$
restricting to a morphism $\PsiDO(\cO_{P'},z'_{P'})_{\leq 0} 
\to \PsiDO(\cO_{P},z_{P})_{\leq 0} $, 
and such that the diagram 
$$
\begin{array}{ccc}
\PsiDO(\cK_{P'},z'_{P'})_{\leq 0} & @> \al >> & \PsiDO(\cK_{P},z_{P})_{\leq 0} 
\\ 
\uparrow{L_{P'}^{z'_{P'}}} & \; & \uparrow{L_{P}^{z_{P}}} 
\\ 
B_{X'}^{C',\rat} & @>>> & B_{X}^{C,\rat}
\end{array}
$$
commutes. The Laurent expansion morphisms behave with respect to 
changes of the vector fields according to (\ref{vfields:Laurent}). 
So we may replace $X'$ by a rational vector field $Y'$, 
without any zero or pole at $P'$. We denote by $Y$ the corresponding 
vector field on $C$, and by $Y_\local,Y'_\local$ the formal expansions of 
$Y,Y'$ at $P,P'$. 

Now we have a commuting diagram 
$$
\begin{array}{ccc}
\PsiDO(\cK_{P'},Y'_\local)_{\leq 0} & @> \al >> & \PsiDO(\cK_{P},Y_\local)_{\leq 0}
\\ 
\uparrow & \; & \uparrow
\\ 
B_{X'}^{C',\rat} & @>>> & B_{X}^{C,\rat}
\end{array}
$$
Since $Y'_\local$ preserves $\cO_{P'}$, $\PsiDO(\cK_{P'},Y'_\local)_{\leq 0}$ 
contains a subalgebra $\PsiDO(\cO_{P'},Y'_{\local})_{\leq 0}$. 
The assumptions on $Y'$ allow to identify $\PsiDO(\cO_{P'},z'_{P'})_{\leq 0}$ 
with $\PsiDO(\cO_{P'},Y'_{\local})_{\leq 0}$. Let us show that 
$\al$ takes this subalgebra to $\PsiDO(\cO_{P'},z'_{P'})_{\leq 0}$. 

The map $\al$ takes $\cO_{P'}$ to $\cO_P$, and it takes $(D_{Y'_\local})^{-1}$ 
to $(D_{Y_\local})^{-1}$. Now we have $Y'_\local = \mu
{\pa\over{\pa z'_{P'}}}$, where  $\mu\in \CC[[z'_{P'}]]^\times$. 
On the other hand, the local expansion of $Y$ at $P$ has the form 
$$
Y_\local = \pi (z_P)^{1 - \nu_P} {\pa\over{\pa z_P}}, 
$$
with $\pi \in \CC[[z_P]]^\times$. This expansion implies that 
$(D_{Y_\local})^{-1}$ has the form 
$$
\sum_{i\geq 1} \pi_i (D_{\pa/\pa z_P})^{-i},
$$
where $\pi_i \in (z_P)^{i + \nu_P -2} \CC[[z_P]]$, so 
$(D_{Y_\local})^{-1}$ belongs to $\PsiDO(\cO_P,z_P)_{\leq 0}$. 

Now $\al$ takes $(D_{\pa/\pa z'_{P'}})^{-1}$ to 
$(D_{Y_\local})^{-1} \vp^*(\mu)$, which belongs to 
$\PsiDO(\cO_P,z_P)_{\leq 0}$. So $\al$ takes the generators of 
$\PsiDO(\cO_P,z_P)_{\leq 0}$ to $\PsiDO(\cO_{P'},z'_{P'})_{\leq 0}$, so 
$$
\al(\PsiDO(\cO_P,z_P)_{\leq 0} ) \subset  \PsiDO(\cO_{P'},z'_{P'})_{\leq 0} . 
$$
This implies that $\vp^*$ takes 
$$
\big( \prod_{P'\in C'} L_{P'}^{z'_{P'}}\big)^{-1} \big( 
\prod_{P'\in C'} \PsiDO(\cO_{P'},z'_{P'})_{\leq 0} \big) 
$$
to 
$$
\big( \prod_{P\in C} L_{P}^{z_{P}}\big)^{-1} \big( 
\prod_{P\in C} \PsiDO(\cO_{P},z_{P})_{\leq 0} \big) ,  
$$
so $\vp^*(B_{X'}^{(C')}) \subset B_{X}^{(C)}$. 

In the same way, one proves that 
$\vp^*(B_{X'}^{(C',D')}) \subset B_{X}^{(C,D)}$. 

This ends the proof of Proposition \ref{prop:funct}. 
\hfill \qed \medskip 

\section{An explicit example: the rational case} \label{sect:example}

It is easy to see that the results of Theorem \ref{thm:pseudo}
also hold in the following cases: $g = 0$, deg$(D) \geq 2$; 
and $g = 1$, deg$(D)\geq 1$. In this section, we study the first case. 

\subsection{Presentation of the classical algebra} 

Let us set $C = \CC P^1$, $D = N\infty$, where $N\geq 2$. 
We have then 
$$
A^{(C,D)} = \oplus_{i\geq 0} A^{(C,D)}_i,
$$ 
where 
$A^{(C,D)}_i = H^0(\CC P^1,K(D)^{\otimes i}) = \{f(z)(dz)^i,$ where 
$f(z)$ is a polynomial of degree $\leq i(N-2)\}$. 

\begin{prop}
$A^{(C,D)}$ may be presented as follows: generators are 
$\omega_a = z^a dz$, $a = 0,\ldots,N-2$, and relations are 
$$
\omega_a \omega_b = \omega_c \omega_d, 
$$
for any quadruple $(a,b,c,d)$ such that $a+b = c+d$. 
\end{prop}

{\em Proof.} Let $\wt A(N)$ be the algebra with generators 
$t_0,\ldots,t_{N-2}$, and relations 
\begin{equation} \label{rel:t}
t_a t_b = t_c t_d, 
\end{equation}
for any quadruple $(a,b,c,d)$ such that $a+b = c+d$. 
Then $\wt A(N)$ is the sum of its homogeneous components $\wt A(N)_i$,  
and relations (\ref{rel:t}) imply that a generating family of 
$\wt A(N)_i$ is given by the union of the 
\begin{equation} \label{gen:fam:t}
(t_0)^\al (t_{N-2})^\beta t_k, \quad \al,\beta\geq,\al+\beta = i-1, 
\quad k = 0,\ldots,N-3, 
\end{equation}
with $t_{N-2}^i$. 

We have an algebra morphism 
\begin{equation} \label{morph:tilde}
\wt A(N) \to A^{(C,D)}, 
\end{equation}
taking each $t_i$ to $\wt \omega_i$. It takes the 
generating family (\ref{gen:fam:t}) to a basis of 
$A_i^{(C,D)}$, which proves as the same time that 
this family is a basis, and that (\ref{morph:tilde}) is an 
isomorphism.  
\hfill \qed\medskip 

The Poisson bracket on $A^{(C,D)}$ is given by 
$\{\omega_a, \omega_b\} = (b-a) z^{a+b-1} (dz)^3$, so in terms of 
generators 
$$
\{\omega_a, \omega_b\} = (b-a) \omega_c \omega_d \omega_e, 
$$
for any $(c,d,e)$ such that $c+d+e = a+b-1$. 

\subsection{Quantized algebras $B_X$ and $B'_X$}

The field of rational functions on $\CC P^1$ is the field of rational 
fractions $\CC(z)$. Equip it with its derivation 
$X = \pa_z = {{d}\over {dz}}$. Then $B_X^{\on{rat}} = 
\PsiDO(\CC(z),\pa_z)_{\leq 0}$. 

Lifts in $B_X$ of the $\omega_a$ are the elements 
$\wt\omega_a = (D_{\omega_a^{-1}})^{-1}$, i.e., 
$$
\wt\omega_a = (\pa_z)^{-1} z^a, \quad a = 0,\ldots,N-2. 
$$
Denote by $B'_X$ the subalgebra of $B_X^{\on{rat}}$ generated by the 
$\wt\omega_a$. Since the $\omega_a$ generate $A^{(C,D)}$, 
$B_X$ is the completion of $B'_X$ with respect to the topology of 
$B_X^{\on{rat}}$. 

\begin{thm} 
For any quadruple $(a,b,c,d)$ such that $0\leq a,b,c,d\leq N-2$, 
$a+b = c+d$ and $b>d$, we have 
\begin{equation} \label{rel:omega:tilde}
\wt \omega_a \wt \omega_b - \wt \omega_c \wt \omega_d = 
(d-b) \wt \omega_a \wt \omega_{b-d} \wt\omega_d. 
\end{equation}
Let us define $C$ as the algebra with generators $t_a, a = 0,\ldots,N-2$
and relations 
$$
t_a t_b - t_c t_d = (d-b) t_a t_{b-d-1} t_d , 
$$
for $a,b,c,d = 0,\ldots,N-2$, such that $a+b = c+d$ and $b>d$. 
Let $I_C$ be the ideal of $C$ generated by the $t_a, a = 0,\ldots,N-2$. 
Set $\wh C = \limm_{\leftarrow n} C / (I_C)^n$. Then there is a 
unique continuous algebra isomorphism 
\begin{equation} \label{morph:pres}
\wh C \to B_X,
\end{equation}
taking each $t_a$ to $\wt \omega_a$. This isomorphism 
induces an algebra isomorphism 
$$
C / \big( \cap_{n\geq 0} (I_C)^n \big)  \to B'_X. 
$$
\end{thm}

{\em Proof.} Let us first prove the relation (\ref{rel:omega:tilde}). 
We have  
$$
\wt \omega_a \wt\omega_b = z^{a+b} (\pa_z + {b\over z})^{-1} (\pa_z)^{-1}, 
$$
so 
\begin{align*}
\wt \omega_a \wt \omega_b - \wt \omega_c \wt \omega_d & = 
z^{a+b} \big( (\pa_z + {b\over z})^{-1} - (\pa_z + {d\over z})^{-1} \big)  (\pa_z)^{-1}
\\ & =  
z^{a+b} (\pa_z + {b\over z})^{-1} {{d-b}\over{z}} (\pa_z + {d\over z})^{-1}  (\pa_z)^{-1}
\\ & 
= (d-b)
z^{a+b-1} (\pa_z + {{b-1}\over z})^{-1} (\pa_z + {d\over z})^{-1}  (\pa_z)^{-1}
\\ & = 
(d-b) z^a (\pa_z)^{-1} z^{b-d-1} (\pa_z)^{-1} z^d (\pa_z)^{-1}
 = (d-b) \wt\omega_a \wt\omega_{b-d-1} \wt\omega_d. 
\end{align*}

We have $\gr(C) = \oplus_{n\geq 0}
(I_C)^n / (I_C)^{n+1}$. We have a morphism of filtered algebras 
$C \to B'_X$. Moreover, we have $\gr(B'_X) = \gr(B_X) = A^{(C,D)}$, so 
we get an algebra morphism $\gr(C) \to \gr(B_X) = A^{(C,D)}$. 

Select in relations (\ref{rel:omega:tilde}), the subset of relations 
corresponding to   $(a,b,a+b,0)$ for $a,b$ such that $a+b \leq N-2$, and 
$(a,b,N-2,a+b-(N-2))$ for $a,b$ such that $a+b >N-2$. Then this subset of relations
implies that a generating family of $\gr_n(B'_X)$ is the union of all 
$$
(\wt\omega_0)^\al (\wt\omega_{N-2})^\beta \wt\omega_i, \quad
i = 0,\ldots,N-3, \; \al + \beta = n-1, 
$$
with $(\wt\omega_{N-2})^i$. 
The morphism $\gr(C) \to A$ takes it to a basis of $A_n$, 
so $\gr(C)\to \gr(B'_X)$ is an isomorphism. This implies that 
the map $\wh C \to B_X$ obtained by completing $C \to B'_X$
is an isomorphism. This fact now implies that $C/ (\cap_{n\geq 0} (I_C)^n)\to B'_X$ is
injective. Since it is obviously surjective, this map is an 
isomorphism. 
\hfill \qed \medskip 

\begin{remark}
We do not know whether $\cap_{n\geq 0} (I_C)^n = 0$, in other words, whether 
$C$ is separated for the topology defined by the powers of $I_C$. 
\end{remark}

\section{Quantization based on Poincar\'e uniformization} 
\label{sect:analytic}

In this section, we assume that $C$ is defined over $\CC$, and that 
we are given a Poincar\'e uniformization of $C$. We denote by 
$\cH$ the Poincar\'e half-plane, and we denote by $\Gamma$ a discrete subgroup of 
$SL_2(\RR)$, such that there is an analytic isomorphism $\cH / \Gamma \to C$. 

We will recall the results of \cite{CMZ} in the Rankin-Cohen brackets (Section 
\ref{sect:CMZ}); we will show how they give rise to a solution $B^\an$ of the 
problem of quantizing the algebra $A$ (Section \ref{sect:sol:anal}),  and 
that this solution is isomorphic to the quantization $B_X$ of Section 
\ref{sect:pseudo} (Section \ref{sect:isom}). For  simplicity, we restrict 
ourselves to the case $D = 0$. 

\subsection{Rankin-Cohen brackets and pseudodifferential operators on 
$\cH$: the results of \cite{CMZ}} \label{sect:CMZ}

Let us denote by $\Hol(\cH)$ the ring of holomorphic functions on the 
Poincar\'e half-plane and by $\tau$ the coordinate on this plane. 
Let us denote by $\pa_\cH$ its derivation 
$d/d\tau$. Consider the algebra $\PsiDO(\Hol(\cH),\pa_\cH)_{\leq 0}$. 
It is a filtered ring, with associated graded 
$$
\bigoplus_{i\geq 0} H^0(\cH,K_\cH^{\otimes i}), 
$$
where $K_\cH$ is the sheaf of differentials on $\cH$. 

$K_\cH$ has a natural section $d\tau$, which induces isomorphisms 
$H^0(\cH,K_\cH^{\otimes i}) \to \Hol(\cH)$. 

The rings $\PsiDO(\Hol(\cH),\pa_\cH)_{\leq 0}$ and $\bigoplus_{i\geq 0}
H^0(\cH,K_\cH^{\otimes i})$ are equipped with natural actions of $SL_2(\RR)$. 
The paper \cite{CMZ} contains the following results: 

\begin{thm} (see \cite{CMZ}) \label{CMZ:1}
There exists a lifting map 
$$
\lift : \bigoplus_{i\geq 0} H^0(\cH,K_\cH^{\otimes i})
\to \PsiDO(\Hol(\cH),\pa_\cH)_{\leq 0}, 
$$
which is $SL_2(\RR)$-equivariant, and expressed by differential operators
acting on $i$-differentials. So the restriction of $\lift$ to 
$H^0(\cH,K_\cH^{\otimes i})$ maps this space to 
$$
\PsiDO(\Hol(\cH),\pa_\cH)_{\leq -i},
$$
and the composed map 
$$
H^0(\cH,K_\cH^{\otimes i}) \to \PsiDO(\Hol(\cH),\pa_\cH)_{\leq -i}
\to \PsiDO(\Hol(\cH),\pa_\cH)_{\leq -i} / 
\PsiDO(\Hol(\cH),\pa_\cH)_{\leq -i-1}
$$ 
is inverse to the natural isomorphism
$\PsiDO(\Hol(\cH),\pa_\cH)_{\leq -i} / 
\PsiDO(\Hol(\cH),\pa_\cH)_{\leq -i-1} \to H^0(\cH,K_\cH^{\otimes i})$. 

If $\omega\in H^0(\cH,K_\cH^{\otimes i})$ has the form $\omega(\tau)(d\tau)^i$, then 
$\lift(\omega)$ has the expression
$$
\lift(\omega) = \omega(\tau) (\pa_\cH)^{-i} + 
\sum_{n>0} \ell_{i,n} \omega^{(n)}(\tau) (\pa_\cH)^{-i-n}, 
$$ 
where $\ell_{i,n}$ are explicit rational numbers. 
\end{thm}

Denote by $\mu$ the product on $\bigoplus_{i\geq 0} H^0(\cH,K_\cH^{\otimes i})$
obtained by transporting the product of $\PsiDO(\Hol(\cH),\pa_\cH)_{\leq 0}$
by the map $\lift$. Since the product on 
$$
\PsiDO(\Hol(\cH),\pa_\cH)_{\leq 0}
$$
is expressed by differential operators, $\mu$ is a $SL_2(\RR)$-invariant 
star-product on $\bigoplus_{i\geq 0} H^0(\cH,K_\cH^{\otimes i})$. 
More precisely, the authors of \cite{CMZ} show: 

\begin{thm} (see \cite{CMZ}) \label{CMZ:2}
Let us denote by $\mu_{ij}^k$ the map 
$H^0(\cH,K_\cK^{\otimes i}) \otimes H^0(\cH,K_\cK^{\otimes j})
\to H^0(\cH,K_\cK^{\otimes i+j+k})$ induced by $\mu$, then the 
$\mu_{ij}^k$ are the Rankin-Cohen brackets: we have 
$$
\mu_{ij}^k(\omega(\tau)(d\tau)^i , \omega'(\tau)(d\tau)^j)
= \sum_{\al,\beta\geq 0} a_{i,j,\al,\beta}^k
\omega^{(\al)}(\tau) \omega^{\prime(\beta)}(\tau) (d\tau)^{i+j+k}, 
$$
for suitable rational numbers $a_{i,j,\al,\beta}^k$. 
\end{thm}

This result immediately implies that the Rankin-Cohen brackets are 
$SL_2(\RR)$-invariant. 

\subsection{Construction of $B^\an$} \label{sect:sol:anal}

Let $C$ be a complex curve, equipped with an isomorphism 
$C\to \cH/ \Gamma$ of analytic manifolds. This isomorphism 
induces an isomorphism 
$$
A \to \bigoplus_{i\geq 0} \big( H^0(\cH,K_\cH^{\otimes i})\big)^\Gamma
$$
of graded algebras and of Poisson algebras. 

\begin{thm}
Set 
$$
B^\an = \big( \PsiDO(\Hol(\cH),\pa_\cH)\big)^\Gamma .  
$$
Then $B^\an$ is a filtered algebra. Its associated Poisson algebra 
is isomorphic to $A$. 
\end{thm}

{\em Proof.} Let us set $(B^\an)^i = \big( \PsiDO(\Hol(\cH),\pa_\cH)_{\leq -i}
\big)^\Gamma$. This obviously defines a grading on $B^\an$, and the image of the 
composed map 
\begin{equation} \label{seq:invt}
(B^\an)^i / (B^\an)^{i+1} \to \PsiDO(\Hol(\cH),\pa_\cH)_{\leq -i}
/ \PsiDO(\Hol(\cH),\pa_\cH)_{\leq -i-1} \to H^0(\cH,K_\cH^{\otimes i})
\end{equation}
is contained in $\big( H^0(\cH,K_\cH^{\otimes i})\big)^\Gamma$. 
So we have a natural map 
\begin{equation} \label{result:invt}
(B^\an)^i / (B^\an)^{i+1} \to \big( H^0(\cH,K_\cH^{\otimes i}) \big)^\Gamma. 
\end{equation}
Since each map of the sequence (\ref{seq:invt}) is injective, 
so is (\ref{result:invt}). It remains to prove that (\ref{result:invt})
is surjective. Denote by $\lift_{i,\Gamma}$ the restriction of $\lift$ 
to $(H^0(\cH,K_\cH^{\otimes i}))^\Gamma$. Then $\lift_{i,\Gamma}$ is a linear 
map 
$$
\big( H^0(\cH,K_\cH^{\otimes i}) \big)^\Gamma \to 
\PsiDO(\Hol(\cH),\pa_\cH)_{\leq -i}. 
$$
According to Theorem \ref{CMZ:1}, the image of $\lift_{i,\Gamma}$
is actually contained in 
$$
(\PsiDO(\Hol(\cH),\pa_\cH)_{\leq -i})^\Gamma = 
(B^\an)^i.
$$ 

The composed map 
$$
\big( H^0(\cH,K_\cH^{\otimes i})\big)^\Gamma 
\stackrel{\lift_{i,\Gamma}}{\to} (B^\an)^i \to (B^\an)^i / (B^\an)^{i+1}
$$
is obviously a section of (\ref{result:invt}), which proves that this map 
is surjective. \hfill \qed \medskip 

The algebra $B^\an$ also has a "star-product" version. 

\begin{prop} \label{prop:star:pdt}
The product $\mu$, defined in terms of Rankin-Cohen brackets (see Theorem 
\ref{CMZ:2}), restricts to a product on 
$$
A = \bigoplus_{i\geq 0} \big( H^0(\cH,K_\cH^{\otimes i})\big)^\Gamma. 
$$ 
The restriction of $\lift$ induces an isomorphism between 
$(A,\mu)$ and $B^\an$.  
\end{prop}

\subsection{Isomorphism with the construction of Section \ref{sect:pseudo}} 
\label{sect:isom}

\begin{prop} \label{prop:isom:2constr}
For any nonzero rational vector field $X$ on $C$, there is an 
isomorphism 
$$
\al_X : B_X \to B^\an
$$
of complete filtered algebras. If $Y$ is another nonzero rational 
vector field on $C$, then $\al_X = \al_Y \circ i_{XY}$.  
\end{prop}

{\em Proof.} Let us denote by $\Mer(\cH)$ the ring of meromorphic functions 
on $\cH$, with only regular singularities. Then $\pa_\cH$
extends to a derivation fo $\Mer(\cH)$ (which we also denote by 
$\pa_\cH$). We set 
$$
\PsiDO(\cH)_{\leq 0}^\mer = \PsiDO(\Mer(\cH),\pa_\cH)_{\leq 0}. 
$$
We will also set 
$$
\PsiDO(\cH)_{\leq 0} = \PsiDO(\Hol(\cH),\pa_\cH)_{\leq 0}. 
$$
Then we have a commuting square of algebras
$$
\begin{array}{ccc}
\big( \PsiDO(\cH)^\mer_{\leq 0}\big)^\Gamma
& \hookrightarrow & \PsiDO(\cH)^\mer_{\leq 0}
\\ 
\uparrow& \; & \uparrow
\\ 
\big( \PsiDO(\cH)_{\leq 0}\big)^\Gamma
& \hookrightarrow & \PsiDO(\cH)_{\leq 0}
\end{array}
$$
where the vertical arrows are injective. We have a natural 
injection $\CC(C) \hookrightarrow \Mer(\cH)$, induced by the projection 
$\cH \to C$. Moreover, let $X$ be a nonzero rational vector field on $C$. 
The lift of $X$ to $\cH$ may be expressed in the form $X(\tau) {d\over{d\tau}}$, 
where $X(\tau) \in \Mer(\cH)$, so Lemma \ref{lemma:pdo} implies that there is a 
canonical morphism 
$$
\PsiDO(\CC(C),X)_{\leq 0} \to \PsiDO(\cH)^\mer_{\leq 0}. 
$$
Lemma \ref{lemma:pdo} also shows that the image of this morphism is 
contained in $\big( \PsiDO(\cH)^\mer_{\leq 0}\big)^\Gamma$. 
Recall that we have $\PsiDO(\CC(C),X)_{\leq 0} = B_X^\rat$; so we have
constructed an algebra morphism $B_X^\rat \to 
\big( \PsiDO(\cH)^\mer_{\leq 0}\big)^\Gamma$. 

We now want to prove that we have a commuting square of algebras
$$
\begin{array}{ccc}
B_X^\rat & \to & 
\big( \PsiDO(\cH)^\mer_{\leq 0}\big)^\Gamma
\\ 
\uparrow& \; & \uparrow
\\ 
B_X & \hookrightarrow & 
\big( \PsiDO(\cH)_{\leq 0}\big)^\Gamma
\end{array}
$$
where the vertical arrows are injective. 
We proceed as follows: 

(a) for any point $P\in C$, there are natural 
Laurent expansion morphisms
$$
\wt L_P^{z_P} : \big( \PsiDO(\cH)^\mer_{\leq 0}\big)^\Gamma 
\to \PsiDO(\cK_P,z_P)_{\leq 0}, 
$$
such that the diagram 
$$
\begin{array}{ccc}
\; & \; & \PsiDO(\cK_P,z_P)_{\leq 0} \\ 
\; & \stackrel{L_P^{z_P}}{\nearrow} & 
\uparrow{\wt L_P^{z_P}}
\\ 
B_X^\rat & \to  & 
\big( \PsiDO(\cH)^\mer_{\leq 0}\big)^\Gamma
\end{array}
$$
commutes.  

(b) $\big( \PsiDO(\cH)_{\leq 0}\big)^\Gamma$ may be characterized as 
the preimage of $\prod_{P\in C} \PsiDO(\cO_P,z_P)_{\leq 0}$ by 
$$
\prod_{P\in C} \wt L_P^{z_P} : 
\big( \PsiDO(\cH)^\mer_{\leq 0}\big)^\Gamma \to 
\prod_{P\in C} \PsiDO(\cK_P,z_P)_{\leq 0}. 
$$

(c) For each $P\in C$, the composed maps 
$$
B_X \to B_X^\rat \stackrel{L_P^{z_P}}{\to} \PsiDO(\cK_P,z_P)_{\leq 0}
$$ 
and
$$
B_X \to B_X^\rat \to \big( \PsiDO(\cH)^\mer_{\leq 0}\big)^\Gamma
\stackrel{\wt L_P^{z_P}}{\to} \PsiDO(\cK_P,z_P)_{\leq 0}
$$ 
coincide, so the image of the latter map is contained in 
$\PsiDO(\cO_P,z_P)_{\leq 0}$. So the image of the composed map  
$$
B_X \to B_X^\rat \to \big( \PsiDO(\cH)^\mer_{\leq 0}\big)^\Gamma
$$ 
is contained in $\big( \PsiDO(\cH)_{\leq 0}\big)^\Gamma$. 
So we have constructed a morphism
\begin{equation} \label{star}
B_X \to  \big( \PsiDO(\cH)_{\leq 0}\big)^\Gamma
\end{equation}
of filtered algebras. 

(d) Both algebras $B_X$ and $\big( \PsiDO(\cH)_{\leq 0}\big)^\Gamma$
are complete and separated for their filtrations. Their associated graded algebras
are isomorphic. Therefore (\ref{star}) induces an algebra isomorphism. 
This proves Proposition \ref{prop:isom:2constr}. \hfill \qed \medskip 

\begin{remark}
The authors of \cite{CMZ} actually define a family of star-products, 
depending on a parameter $\kappa$. In the language of Section \ref{line}, 
this construction corresponds to replacing the algebra $B_X^{(C,D)}$
by the family of algebras $B_X^{(C,D),\la}$, where the generalized line 
bundle $\la$ is $\kappa \bar\al$, and $\bar\al$ is an element of $\CC C$
such that the class of $2\bar\al$ modulo $\res(\CC(C)^\times)$ is equal to
the canonical bundle $K_C$. 
\medskip \end{remark}

\begin{remark}
To be able to use Proposition \ref{prop:star:pdt}, one needs to 
know the group $\Gamma$ corresponding to a given curve $C$. 
This is the case, by definition, if $C$ is a modular curve. 
In this case, a classical problem is to find algebraic equations 
for this curve. This problem is solved using the algebra of modular 
forms. The corresponding "quantum" problem is to give a presentation of the 
algebra $B_X$ (or equivalently, of $\oplus_{i\geq 0} \big( 
H^0(\cH,K_\cH^{\otimes i}) \big)^\Gamma$, equipped with its Rankin-Cohen 
star-product structure $\mu$). 
\end{remark}

\section{Differential liftings} \label{sect:lifting}

The lifting 
$$
\bigoplus_{i\geq 0} H^0(C,K^{\otimes i}) \to B_X
$$
constructed in the proof of Theorem \ref{thm:pseudo} (see step 
(c) of Section \ref{sect:proof}) relies on estimation of the dimensions 
of cohomology groups. Contrary to the operation $\lift$ of Theorem \ref{CMZ:1}, it is 
therefore not a local operator. We now study the problem of constructing 
such a local, or differential, lifting, in the algebraic framework.
We will prove that the set $\Lift_\diff(C)$ of such liftings is a torsor under the 
action of a group $\Aut_\diff(C)$. Poincar\'e uniformization 
yields a point of this torsor. We do not know an algebraic 
way to construct a point of the torsor $\Lift_\diff(C)$, but we study
some algebraic structures provided by such a point.  

\subsection{Differential liftings}

A differential lifting of the isomorphism $\gr(B_X) \to A$
is defined as the following data: for each rational vector field $X$, 
this is a collection $(\La^X_{i,j})_{i,j\geq 0}$ of rational differential 
operators $\La^X_{i,j} : \{$rational $i$-differentials$\}\to \CC(C)$. 
This collection is subject to the following conditions: 

(1) define $\La_X : A^\rat \to B_X^\rat$ by 
$$
\underline{\omega} = (\omega_i)_{i\geq 0}
\mapsto \La^X(\underline{\omega}) = \sum_{i,j\geq 0} 
\La^X_{i,j}(\omega_i) (D_X)^{-j}, 
$$
then $\La^X$ maps $A^\rat_i$ to $(B_X^\rat)^i$, the composed map 
$A^\rat_i \stackrel{\La^X}{\to} (B_X^\rat)^i \to (B_X^\rat)^i / (B_X^\rat)^{i+1}$
is the inverse of the canonical map, and $\La^X(1) = 1$;  

(2) for any pair $X,Y$ of nonzero vector fields, we have 
$i_{X,Y}\circ \La^X = \La^Y$

(3) condition (2) implies that for any $P\in C$, 
$\La^X$ induces a map  
$$
\La^X_P : \bigoplus_{i\geq 0} \CC((z_P))(dz_P)^i
\to \PsiDO(\cK_P,z_P)_{\leq 0}. 
$$
Then for any $P\in C$, $\La^X_P$ maps $\bigoplus_{i\geq 0} \CC[[z_P]](dz_P)^i$
to $\PsiDO(\cO_P,z_P)_{\leq 0}$. 

(If a nonzero vector field $X_0$ is fixed, then for any
family $(\La^{X_0}_{i,j})_{i,j}$ satisfying conditions 
(1), (3) for $X_0$, condition (2) uniquely determines 
a differential lifting $(\La^X_{i,j})_{i,j\geq 0}$ extending 
 $(\La^{X_0}_{i,j})_{i,j}$.)

Conditions (1), (2) and (3) imply immediately that $\La^X$ induces 
a linear map 
$$
\rho(\La^X) : A \to B_X, 
$$
which is a section of the canonical map $\gr(B_X) \to A$, and therefore
induces an isomorphism $\wh{\rho(\La^X)} : \wh A \to B_X$, where
$\wh A$ is the completion $\prod_{i\geq 0} A_i$. 

Let us denote by $\Lift_\diff(C)$ the set of all differential 
lifts on $C$. For any nonzero rational vector field $X$, 
the assignment $\La^X \mapsto \wh{\rho(\La^X)}$ is a map 
$$
\rho : \Lift_\diff(C) \to \on{Isom}(\wh A,B_X).
$$
We will now see that both sides of this map are 
principal homogneous spaces (torsors) and that $\rho$
is a morphism of torsors.

\subsection{The group $\Aut_\diff(C)$}

Define $\DO(K^{\otimes i},K^{\otimes j})$ as the space of all regular 
differential operators on $C$, from $K^{\otimes i}$ to $K^{\otimes j}$. 
Define $\DO(K^{\otimes i},K^{\otimes j})_{\leq k}$ as the subspace of 
all such operators of order $\leq k$. Set 
$$
\gr_k \big( \DO(K^{\otimes i},K^{\otimes j}) \big) = 
\DO(K^{\otimes i},K^{\otimes j})_{\leq k} / 
\DO(K^{\otimes i},K^{\otimes j})_{\leq k-1}. 
$$
Then we have a graded linear injection 
\begin{equation} \label{inj:DO}
\bigoplus_{k\geq 0} \gr_k\big( \DO(K^{\otimes i},K^{\otimes j}) \big) 
\hookrightarrow  
\bigoplus_{k\geq 0} H^0(C,K^{\otimes j-i-k}) . 
\end{equation}
It follows that when $i>j$, $\DO(K^{\otimes i},K^{\otimes j}) = 0$, 
and if $i = j$,  $\DO(K^{\otimes i},K^{\otimes j}) = \CC$. 

Define $\End_{\diff}(C)$ as follows 
$$
\End_\diff(C) = \wh \bigoplus_{i\leq j} \DO(K^{\otimes i},K^{\otimes j}) ,   
$$
where $\wh\oplus$ is the completed direct sum (direct product). 
Then composition of differential operators induces an algebra structure on 
$\End_{\diff}(C)$. Projection of the diagonal summands induces an 
algebra morphism $\End_\diff(C) \to \prod_{i\geq 0} \CC$. The preimage of 
$\prod_{i\geq 0} 1$ in $\End_\diff(C)$ is a group, which we
denote $\Aut_\diff(C)$. It is easy to see that this is a prounipotent 
algebraic group, as is the subgroup $\Aut_{\diff,1}(C)$ of elements preserving $1$. 

Define $\Aut(\wh A)$ as the group of all continuous linear automorphisms
of $\wh A = \wh\oplus_{i\geq 0} A_i$.

\begin{prop} 
There is a natural group morphism $\Aut_{\diff,1}(C) \to \Aut(\wh A)$. 
The map $\rho$ is a torsor morphism, compatible with this group 
morphism. 
\end{prop}

We have already mentioned that Poincar\'e uniformization 
provides an element of $\Lift_\diff(C)$. On the other 
hand, $\Lift_\diff(C)$ is a purely algebraic object, so one would 
like an algebraic construction of its elements. We will not give 
such a construction, but only indicate that such elements give rise to 
affine spaces over spaces of differential operators
(Section \ref{sect:lifts?}).

\subsection{} \label{sect:lifts?}

Let us now describe the possible form of a differential lift $\La$. 
If $\beta$ is a rational differential, we may set 
$$
\La(\beta) = (D_{\beta^{-1}})^{-1} ,  
$$
in the notation of Section \ref{rem:indep}. In other words, 
if $\al_0$ is a nonzero rational differential, and $X_0$ is the vector 
field inverse to $\al_0$, we have 
$$
\La(\beta) = (D_{X_0})^{-1} f, 
$$
where $f = \beta / \al_0$. 

Let now $\beta$ be a rational quadratic differential. 
Set $f = \beta / (\al_0)^2$. The element 
$$
(D_{X_0})^{-2} f + {1\over 2} (D_{X_0})^{-3} X_0(f)
$$
of $(B^\rat)^2 / (B^\rat)^4$ is independent on the choice of 
$\al_0$. On the other hand, one can show that there is no expression
$P(\al_0,f)$ of the form $\sum_n a_n (X_0)^n(f)$, such that the element 
$$
\La(\beta) = (D_{X_0})^{-2} f + {1\over 2} (D_{X_0})^{-3} X_0(f)
+ (D_{X_0})^{-4} P(\al_0,f)
$$
of $(B^\rat)^2 / (B^\rat)^5$ is independent on the choice of 
$\al_0$. So the determination of the coefficient $P(\al_0,f)$
depends on additional data. The space of all possible expressions 
$P(\al_0,f)$ is an affine space, with associated vector 
space $\DO(K^{\otimes 2}, K^{\otimes 4})$. This structure of 
affine space may be viewed as a part of the torsor structure of 
$\Lift_\diff$.

\begin{remark} {\it On the size of $\DO(K^{\otimes n},K^{\otimes m})$}
\label{rem:size}
The injection (\ref{inj:DO}) is not alway surjective: for 
example, where $n=1,m=2,k=1$, the preimage of $1\in H^0(C,\cO_C)$
is the class of all regular connections on $K$; but there is no 
such connection. Indeed, a rational connection on $K$ has the form 
$$
\omega \mapsto \nabla^{\al,\la} (\omega) = \al d(\omega/\al) + \la\omega, 
$$
for $\al,\la$ rational differentials, $\al\neq 0$. The condition for 
$\nabla^{\al,\la}$ to be regular is that $\la$ has poles of order $\leq 1$
at each $P\in C$, and if $(\al) = \sum_{P\in C} n_P P$, then $\res_P(\la) 
= - n_P$. Such a $\la$ cannot exist, by the residue theorem (which says 
that $\sum_{P\in C} \res_P(\la) = 0$, whereas $\sum_{P\in C} n_P = 2g-2$).  
\end{remark}

\section{Concluding remarks} \label{sect:final}

\subsection{The elliptic case}

When $g = 1$ and the degree of $D$ is $>0$, the above construction 
of the algebra $B_X^{(C,D)}$ may still be carried out. Its classical 
limit is the algebra $A^{(C,D)}$. Let us compare them to the elliptic 
algebras of \cite{FO}. 

$A^{(C,D)}$ may be described as follows: $A^{(C,D)} = \oplus_{i\geq 0}
A_i^{(C,D)}$, where $A_i^{(C,D)} = H^0(C,\cO(D)^{\otimes i})$. 
We view elements of $A_i^{(C,D)}$ as rational functions on $C$, 
with divisor $\geq -i D$. In particular, the derivation $f\mapsto {d\over{dz}}(f) = f'$
can be applied to these functions (here $z$ is a uniformizing parameter of 
$C$). The algebra structure of $A^{(C,D)}$ is graded and induced
by the product of rational functions. Its Poisson bracket is defined
as follows: it is homogeneous of degree $1$, and if $f\in A_i^{(C,D)}$, 
$g\in A_j^{(C,D)}$, then $\{f,g\}\in A_{i+j+1}^{(C,D)}$ corresponds to 
\begin{equation} \label{PB:funs}
i (f'g)(z) - j (fg')(z).  
\end{equation}  
It turns out that for any integer $d\geq 0$, one can define a 
Poisson structure on the algebra $A^{(C,D)}$, by requiring that 
it is homogeneous of degree $d$, and for $f\in A_i^{(C,D)}$, 
$g\in A_j^{(C,D)}$,  $\{f,g\}\in A_{i+j+d}^{(C,D)}$ corresponds to 
(\ref{PB:funs}).  The structure studied in this paper corresponds to 
$d = 1$, and the structure of \cite{FO} corresponds to $d=0$. 

As we have said, the quantization of the first structure may be 
done in terms of pseudodifferential operators. The quantization $A_{\FO}$ 
of the second structure was achieved in \cite{FO}. It can be expressed 
in terms of difference operators: if $\sigma = e^{\al (d/dz)}$ 
is a translation of $C$, elements of $(A_{\FO})_n$ are operators
of the form $f e^{n \al (d/dz)}$, where $f$ is a section of 
$\cO(D+ \sigma(D) +\cdots + \sigma^{n-1}(D))$. $A_{\FO}$ is then 
a graded algebra. 

We do not know a quantization of the Poisson algebras corresponding to 
other values of $d$.

\subsection{Higher-dimensional Poisson structures}

Let us set $A_1 = H^0(C,K)$, then we have a map 
\begin{equation} \label{factor:map}
S^\bullet(A_1) \to A = \bigoplus_{i\geq 0} H^0(C,K^{\otimes i}). 
\end{equation}
When $g = 3,4,5$, one can define a Poisson structure on the algebra
$S^\bullet(A_1)$, such that (\ref{factor:map}) is Poisson. In that case, 
the quantization of $S^\bullet(A_1)$ and of the morphism (\ref{factor:map})
is not known. 

In the other cases, a Poisson structure on $S^\bullet(A_1)$, 
such that (\ref{factor:map}) is Poisson, is not known. One $2$-dimensional 
symplectic leaf of such a Poisson structure would be given by the dual to 
the map (\ref{factor:map}), so it would be isomorphic to the cone
$\Cone(C,D)$. One could try to construct geometrically higher 
dimensional symplectic leaves of this Poisson structure by 
first understanding their geometric interpretation when 
$g=3,4,5$. 

\subsection{Relation with Kontsevich quantization}

When $g = 3,4,5$, one may apply Kontsevich quantization to the algebras
$S^\bullet(A_1)$. Under this quantization, the Poisson central elements
$Q_1,\ldots,Q_{g-2}$ are deformed to central elements. So factoring 
them out gives rise to a quantization $S^\bullet(A_1)_\hbar \to A_\hbar$
of the map (\ref{factor:map}). It is natural to expect that $A_\hbar$ and $B_X$
are isomorphic.  

\subsection{Relation with the Beauville hamiltonians}

In \cite{Beauville}, Beauville introduced integrable systems on 
symmetric powers of $K3$ surfaces. An analogous construction is 
the following. Let $k$ be an integer, $(C,D)$ be the pair of a curve and 
an effective divisor, and $\omega_1,\ldots,\omega_k$ be elements of 
$A^{(C,D)}_1$. Set 
$$
A^{(k)} = \big( (A^{(C,D)})^{\otimes k}\big)^{\SG_k}. 
$$
Then $A^{(k)}$ is a Poisson algebra. For $\phi\in A^{(C,D)}$, 
denote by $\phi^{(i)}$ be image of $\phi$ in the $i$th copy of 
$(A^{(C,D)})^{\otimes k}$. Denote by $\psi_0,\ldots,\psi_k$
the minors of the matrix 
$$
\pmatrix 
\omega_1^{(1)} & \cdots & \omega_k^{(1)} &   1     \\
   \vdots      &   \;   &     \vdots     & \vdots \\ 
\omega_1^{(k)} & \cdots & \omega_k^{(k)} &   1 
\endpmatrix 
$$
Set $H_i = \psi_i / \psi_0$, for $i = 1,\ldots,k$. Then the $H_i$
are a Poisson-commuting family of elements of $\on{Frac}(A^{(k)})$. 
It would be interesting to study the quantization of this 
family using the algebras of pseudodifferential operators introduced
here.

\subsection*{Acknowledgements} 

{\small
We would like to thank V. Rubtsov and M. Olshanetsky for discussions
on the subject of this work. We would also like to thank the MPIM-Bonn, 
as well as IRMA (Strasbourg, CNRS), for support at the time this work 
was being done. The work of A.O. is also partially supported by grants 
RFBF 99-01-01169, RFBR 00-15-96579, CRDF RP1-2254 and INTAS-00-00055. 
}

\end{document}